%%%%%%%%%%%%%%%%%%%%%%%%%%%%%%%%%%%%%%%%%%%%%%%%%%
%%                                              %%
%%                                              %%
%%  "Quantum cyclotomic orders of 3-manifolds"  %%
%%        by Tim Cochran and Paul Melvin        %%
%%                                              %%
%%                                              %%
%%%%%%%%%%%%%%%%%%%%%%%%%%%%%%%%%%%%%%%%%%%%%%%%%%

\input amstex.tex
\input amsppt.sty
\input epsf.def.bmc

%%%%%%% paper format %%%%%%%
\magnification 1200
\nologo
\NoBlackBoxes
\pagewidth{5.2 in}
\pageheight{7.3 in}
\hcorrection{.12 in}

%%%%%%% style macros (sections, theorems, etc.) %%%%%%%
\def\section#1{\penalty-1000\bigskip\noindent{\sl #1}\smallskip\nobreak}
\def\subsection#1{\penalty-1000\medskip\noindent{\sl #1}\smallskip\nobreak}
\def\theorem#1{\medskip \bf Theorem #1.\it} 
\def\lemma#1{\medskip \bf Lemma #1.\it} 
\def\proposition#1{\medskip \bf Proposition #1.\it}
 
\def\proof{{\smallskip\sl Proof. \rm}}
\def\defn#1{\medskip \bf Definition #1.\rm} 
\def\remark#1{\medskip \bf Remark #1.\rm}
\def\question#1{\medskip \bf Question #1.\it}

\def\remarks#1{\medskip \bf Remarks #1.\rm}

\def\claim{\medskip \bf Claim. \ \rm} 
\def\pf#1{{\smallskip\sl Proof of \rm(#1).}} 
\def\endtheorem{\rm\smallskip}
\def\endlemma{\rm\smallskip}
\def\endproposition{\rm\smallskip}

\def\endproof{\qed\medskip}
\def\enddefn{\medskip}
\def\endremark{\medskip}
\def\endremarks{\medskip}

\def\endquestion{\rm\smallskip}

\def\endclaim{\rm\smallskip}
\def\endpf{\qed\medskip}

%%%%%%% figures for sections 4 and 5 %%%%%%%
\def\put#1#2#3{\text{\kern-#1pt{\raise#2pt\hbox{\rlap{#3}}}\kern#1pt}}
\def\dy#1{\text{\kern-4pt\lower5pt\hbox{\epsfscale #1 \epsfbox{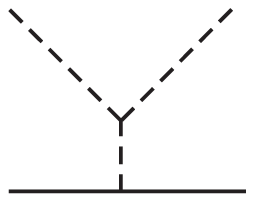}}}}
\def\dloop#1{\text{\kern-4pt\lower5pt\hbox{\epsfscale #1 \epsfbox{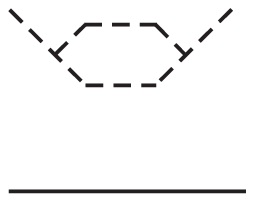}}}}
\def\un#1{\text{\kern-4pt\lower13pt\hbox{\epsfscale #1 \epsfbox{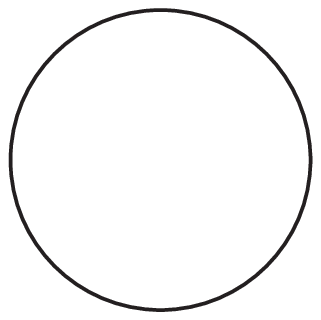}}}}
\def\unt#1{\text{\kern-4pt\lower5pt\hbox{\epsfscale #1 \epsfbox{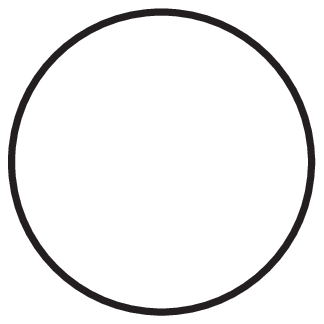}}}}
\def\uni#1{\text{\kern-4pt\lower4pt\hbox{\epsfscale #1 \epsfbox{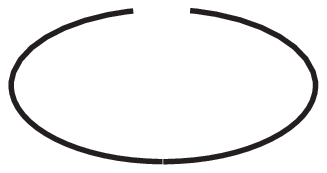}}}}
\def\unii#1{\text{\kern-4pt\lower5pt\hbox{\epsfscale #1 \epsfbox{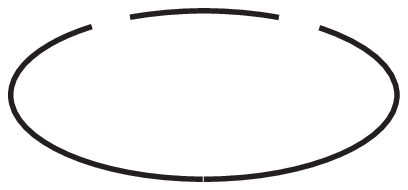}}}}
\def\hopf#1{\text{\kern-4pt\lower5pt\hbox{\epsfscale #1 \epsfbox{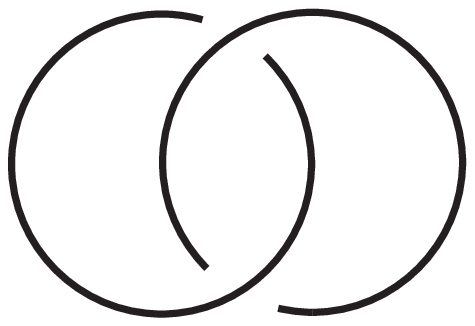}}}}
\def\dhopf#1{\text{\kern-4pt\lower5pt\hbox{\epsfscale #1 \epsfbox{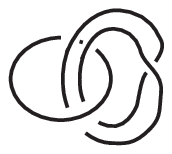}}}}
\def\bo#1{\text{\kern-4pt\lower13pt\hbox{\epsfscale #1\epsfbox{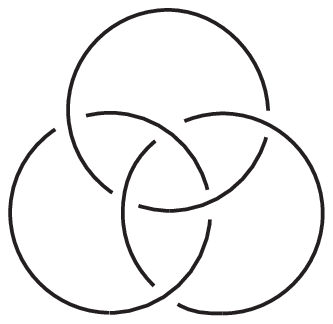}}}}
\def\tri#1{\text{\kern-4pt\lower13pt\hbox{\epsfscale #1 \epsfbox{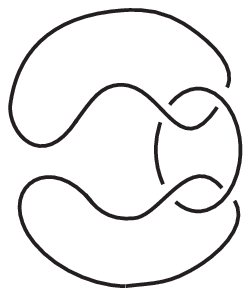}}}}
\def\dtri#1{\text{\kern-4pt\lower13pt\hbox{\epsfscale #1 \epsfbox{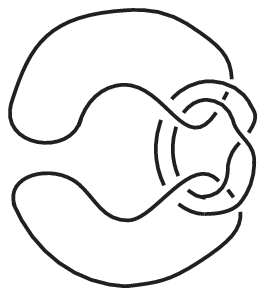}}}}
\def\wh#1{\text{\kern-4pt\lower13pt\hbox{\epsfscale #1 \epsfbox{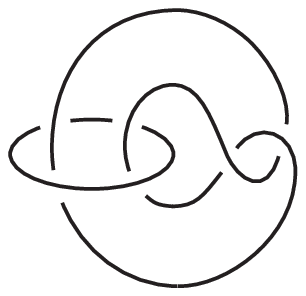}}}}
\def\hop#1{\text{\kern-4pt\lower13pt\hbox{\epsfscale #1 \epsfbox{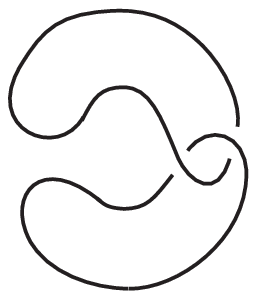}}}}
\def\boa#1{\text{\kern-4pt\lower13pt\hbox{\epsfscale #1 \epsfbox{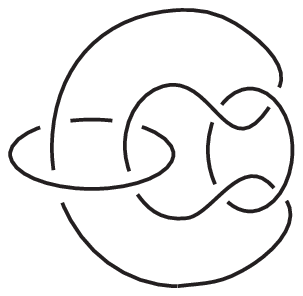}}}}
\def\dboa#1{\text{\kern-4pt\lower13pt\hbox{\epsfscale #1 \epsfbox{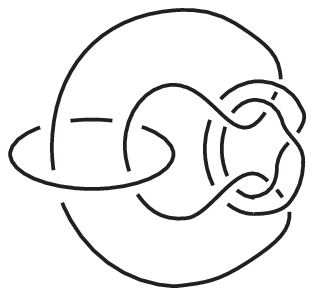}}}}
\def\dboat#1{\text{\kern-4pt\lower13pt\hbox{\epsfscale #1
    \epsfbox{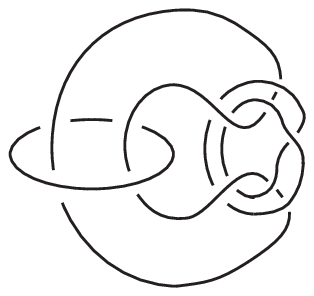}}}}
\def\caps#1{\text{\kern-4pt\lower13pt\hbox{\epsfscale #1 \epsfbox{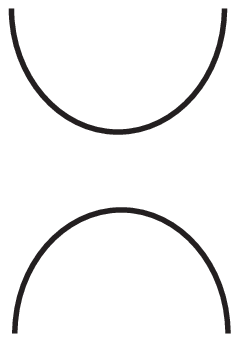}}}}
\def\cross#1{\text{\kern-4pt\lower5pt\hbox{\epsfscale #1 \epsfbox{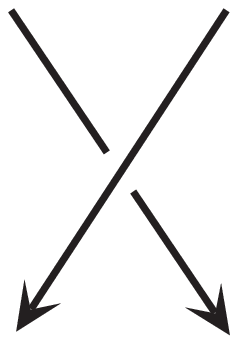}}}}
\def\down#1{\text{\kern-4pt\lower5pt\hbox{\epsfscale #1 \epsfbox{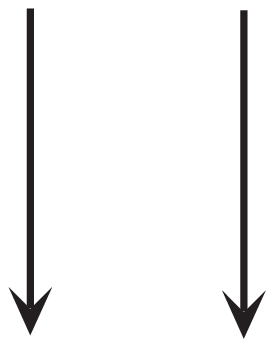}}}}
\def\orcaps#1{\text{\kern-4pt\lower5pt\hbox{\epsfscale #1 \epsfbox{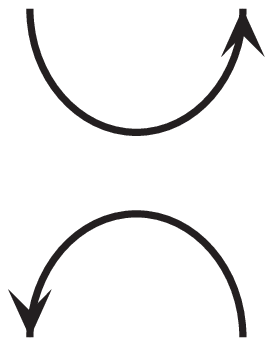}}}}
\def\spi#1{\text{\kern-4pt\lower13pt\hbox{\epsfscale #1 \epsfbox{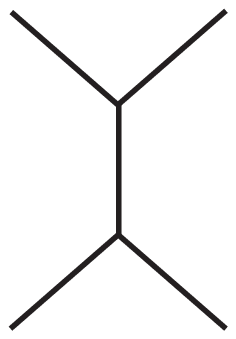}}}}
\def\spis#1{\text{\kern-4pt\lower14pt\hbox{\epsfscale #1 \epsfbox{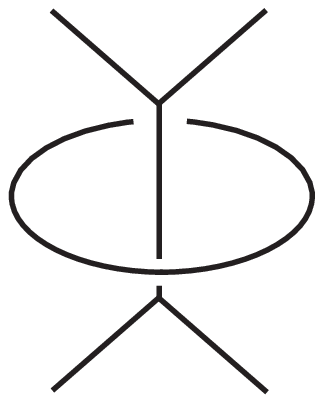}}}}
\def\pos#1{\text{\kern-4pt\lower2pt\hbox{\epsfscale #1 \epsfbox{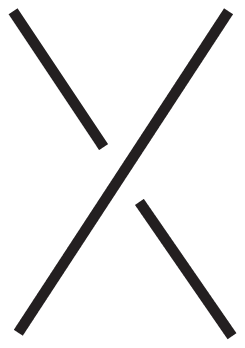}}}}
\def\neg#1{\text{\kern-4pt\lower2pt\hbox{\epsfscale #1 \epsfbox{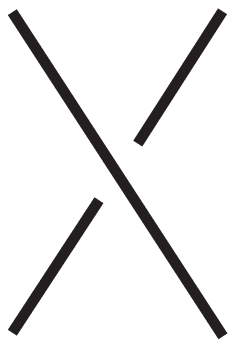}}}}
\def\inj#1{\text{\kern-4pt\lower2pt\hbox{\epsfscale #1 \epsfbox{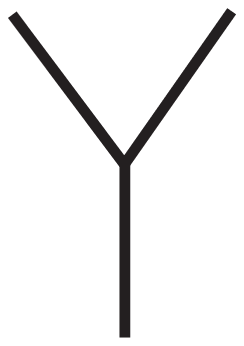}}}}
\def\pro#1{\text{\kern-4pt\lower2pt\hbox{\epsfscale #1 \epsfbox{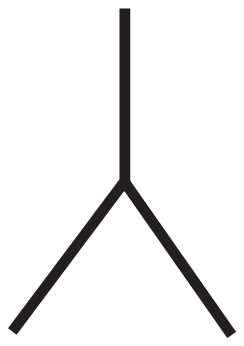}}}}
\def\ring#1{\text{\kern-4pt\lower13pt\hbox{\epsfscale #1 \epsfbox{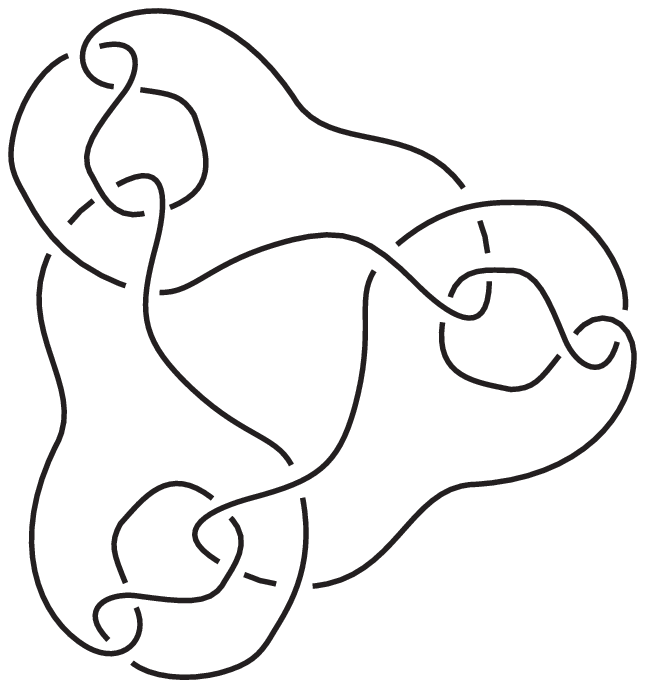}}}}
\def\strand#1{\text{$\overset#1\to{\put{0}{0}{$|$}
  \put{0}{8}{$|$} \put{0}{-8}{$|$}}$}}
\def\strandi#1{\text{$\overset#1\to{\put{0}{0}{$|$}
  \put{0}{9}{$|$} \put{0}{-13}{$|$}}$}}

%%%%%%% standard macros %%%%%%%
\def\Z{\Bbb Z}
\def\Q{\Bbb Q}

\def\C{\Bbb C}
\def\min{\text{\rm min}}
\def\max{\text{\rm max}}

\def\and{\text{\rm \ and }}

\def\rk{\text{\rm rk}}
\def\iff{\Longleftrightarrow}

%%%%%%% special macros %%%%%%%
\def\br#1{\langle#1\rangle}
\def\approxf{\underset f\to\approx}
\def\A{\Cal A}
\def\M{\Cal M}
\def\P{\Bbb P}
\def\Zp{\Z_p}
\def\Zpk{\Z_{p^k}}
\def\po{P_{\text{odd}}}
\def\ppo{\P_{\text{odd}}}
\def\G{\Gamma}
\def\L{{\Lambda}}
\def\lp{{\L_p}}
\def\qp{{Q_p}}
\def\ak{(a,k]}
\def\ac{(a|c)}
\def\kc{[k,c)}
\def\oo{{\text{\rm o}}}
\def\o{{\frak o}}

\def\op{{\frak o_{p}}}
\def\osp{{\frak o_{\!/\!p}}}
\def\b{{b}}
\def\bp{{b_p}}
\def\rp{{r_p}}
\def\sp{{s_p}}
\def\lk{\text{\sl lk}}
\def\v{{\text{v}}}
\def\vp{{\v_p}}
\def\calo{{\Cal O}}
\def\calodv{{\calo^d_\v}}
\def\calodp{{\calo^d_{\!p}}}
\def\calodsp{{\calo^d_{\!/\!p}}}
\def\pc{{\varphi_p}} %% pth cyclotomic polynomial %% 

%%%%%%% references %%%%%%%
\def\cite#1{[#1]}
\def\BHMV{BHMV} % blanchet-..., remarks on the ... % 
\def\BS{BS} % borevich-shafarevich no thy % 
\def\Co{Co} % Cochran homotopy theory and link conc inv %
\def\CGO{CGO} % cochran-gerges-orr % 
\def\CMI{CM1} % cochran-melvin finite type % 
\def\CMII{CM2} % cochran-melvin milnor degree % 
\def\Dr{Dr} % drinfeld quasi hopf algebras, leningrad % 
\def\Ga{Ga} % garoufalidis thesis %
\def\Ge{Ge} % gerges thesis %
\def\HM{HM} % habegger-masbaum %
\def\IR{IR} % ireland-rosen no thy %
\def\Je{Je} % jeffrey lens %
\def\Jo{Jo} % jones poly %
\def\Ka{Ka} % kauffman bracket %
\def\Ki{Ki} % kirby calculus % 
\def\KL{KL} % kauffman-lins recoupling book % 
\def\KMI{KM1} % kirby-melvin inventiones % 
\def\KMII{KM2} % kirby-melvin lens % 
\def\KR{KR} % kirillov-reshetikhin, ref Kac ed World Sci volume % 
\def\KS{KS} % kricker-spence finite type % 
\def\La{La} % lawrence asymptotic exp %
\def\LI{L1} % lickorish first appearance of omega %
\def\LII{L2} % lickorish skein method %
\def\Le{Le} % le universal %
\def\LMI{LM1} % le-murakami homfly %
\def\LMII{LM2} % le-murakami parallel %
\def\LMO{LMO} % le-murakami-ohtsuki %
\def\LW{LW} % lin-wang ohtsuki %
\def\MR{MR} % masbaum-roberts %
\def\MW{MW} % masbaum-wenzl %
\def\Mi{Mi} % milnor isotopy of links %
\def\MI{M1} % murakami su(2) %
\def\MII{M2} % murakami so(3) %
\def\OI{O1} % ohtsuki su(2) series %
\def\OII{O2} % ohtsuki so(3) series %
\def\RTI{RT1} % reshetikhin-turaev ribbon % 
\def\RTII{RT2} % reshetikhin-turaev inv % 
\def\RI{R1} % rozansky reshet form for jones and milnor inv % 
\def\RII{R2} % rozansky p-adic convergence % 
\def\Ta{Ta} % takata su(n) % 
\def\Tu{Tu} % turaev's book % 
\def\Wa{Wa} % washington cyclotomic %
\def\Wi{Wi} % witten jones %

\topmatter

\title Quantum cyclotomic orders of $3$-manifolds \endtitle

\author 
Tim D. Cochran
and 
Paul Melvin
\endauthor

\date September,
 1998
\enddate

\address{Mathematics Department, Rice University, Houston, TX 77251, U.S.A.}
\endaddress

\email{cochran\@rice.edu}
\endemail

\address{Mathematics Department, Bryn Mawr College, Bryn Mawr, PA 19010, U.S.A.}
\endaddress

\email{pmelvin\@brynmawr.edu}
\endemail

\thanks
Both authors gratefully acknowledge the support of Research Professorships at the
Mathematical Sciences Research Institute, Berkeley, California.  The first author
was also partially supported by the National Science Foundation.
\endthanks

\abstract{This paper provides a topological interpretation for number theoretic
properties of quantum invariants of $3$-manifolds.  In particular, it is shown that
the $p$-adic valuation of the quantum $SO(3)$-invariant of a $3$-manifold $M$, for odd
primes $p$, is bounded below by a linear function of the mod $p$ first betti number of
$M$.  Sharper bounds using more delicate topological invariants are given as well.}
\endabstract

\comment


\toc
\head \S0. \ Introduction \endhead
\head \S1. \ Laurent polynomials \endhead
\head \S2. \ Relations with the Kontsevich integral \endhead
\head \S3. \ Cyclotomacy \endhead
\head \S4. \ Quantum invariants \endhead
\head \S5. \ Examples \endhead
\endtoc

\endcomment

\endtopmatter

\document
\baselineskip=12 pt plus 1pt minus .5pt

\heading{\S0. Introduction}
\endheading

Since the birth of quantum topology in the last decade \cite\Jo\,\cite\Wi, one of
the fundamental problems facing topologists has been to find topological
interpretations for the vast array of quantum invariants that have come to
light.  One common characteristic among these invariants is their rich number
theoretic content, and it has been a special challenge to understand which
aspects of this number theory have topological significance.  

Of central interest are the quantum invariants of $3$-manifolds that arise from
the representation theory of classical Lie groups \cite\Wi\,\cite\RTII.  Roughly
speaking there is a complex valued invariant $\tau_p^G$ associated with any
suitable Lie group $G$ (the {\it gauge group}) and integer $p$ (the {\it coupling
constant} or {\it level}).  Typically $\tau_p^G$ takes values in a cyclotomic
field determined by $G$ and $p$, and this is where the number theory comes into
play.  There are of course many algebraic invariants in cyclotomic fields.  Among
the most fundamental are the valuations associated with prime ideals in their
rings of integers.  The object of this paper is to demonstrate a connection
between these valuations and basic topological invariants.  

The results obtained here will also be used in the authors' forthcoming paper on
a theory of finite type invariants for arbitrary $3$-manifolds \cite\CMI.  They
provide the basis for the construction of a rich family of new
invariants of finite type.

For simplicity we shall limit our investigations to the case when $G = SO(3)$ and
$p$ is an odd prime.  Similar considerations should apply to the gauge group
$SU(2)$, since the $SU(2)$ and $SO(3)$ invariants are proportional by a factor
involving only classical homotopy theoretic invariants \cite\KMI, and it is
expected that results will soon follow for other Lie groups at prime levels. 
Extending to composite levels may require new ideas.

So consider the quantum $SO(3)$ invariant at an odd prime level $p$, which
will be denoted simply by $\tau_p$.  This invariant (normalized as in \S4)
takes values in the cyclotomic field $\qp$ of $p$th roots of unity.  Let $\lp$
denote the ring of integers in $\qp$.  There is a distinguished valuation on
$\qp$, associated with the unique prime ideal $H$ in $\lp$ containing $p$, and
its value on any $x\in \qp$ will be called the {\it $p$-order} of $x$.  Now define
the {\it quantum $p$-order} $\op(M)$ of a $3$-manifold $M$ to be the $p$-order of
$\tau_p(M) \in \qp$.  Alternatively $\op(M)$ can be defined as the largest integer
$m$ for which $\tau_p(M) \in H^m$; an easy exercise from the definition of
$\tau_p$ shows that such integers exist.

A beautiful theorem of Hitoshi Murakami \cite\MI\,\cite\MII\ and Masbaum-Roberts
\cite\MR\ states that $\op(M)\ge0$, that is $\tau_p$ actually takes values in
$\lp$.  (Similar results hold for many of the other classical Lie groups by recent
work of Takata \cite\Ta\ and Masbaum-Wenzl \cite\MW.)  It was also shown by
Murakami \cite\MII\ that $\op(M)=0$ if and only if $M$ is a $\Zp$-homology sphere.
This was the first indication that topological information might be carried by
$\op(M)$. 

In this paper it will be shown that the quantum $p$-order of $M$ is bounded below
by a linear function of its mod $p$ first betti number. In particular 
$$
\op(M) \ge \bp(M)\,n/3
$$
where $\bp(M)=\rk(H_1(M;\Zp))$ and $n=(p-3)/2$.  Furthermore, it will be seen that
this is the best possible betti number bound for the order when $\bp\equiv0\pmod
3$.   For $\bp = 1$ or $2$, an improved lower bound $\op \ge n$ will be established
and shown to be sharp.  It is also shown how to strenghthen these bounds by
considering a more refined topological invariant, the ``Milnor degree" (see \S2,
Theorem 4.3 and Remark 4.6).   
    
Along the way, a family of $\Zpk$-valued $3$-manifold invariants $\tau_p^d$
will be introduced.  Using the methods of this paper, it can be shown that these
invariants are of ``finite type" and converge to $\tau_p$.  One striking
consequence of this fact is that any two $3$-manifolds with unequal quantum
$SO(3)$ invariants can be distinguished by finite type invariants (see \cite\CMI\
for details, and \cite\RII\ for the special case of rational homology spheres). 

The paper is organized as follows.  In \S1 we describe two families of Laurent
polynomials, quantum integers and Jones polynomials (including
Ohtsuki's version \cite\OI), and use these to define the $p$-bracket of a
framed link.  The $p$-bracket is the key ingredient in the definition of the
quantum invariant $\tau_p$, given in \S4 along with the definition of the
invariants $\tau_p^d$.  The relationship of the Jones and Ohtsuki polynomials with
the Kontsevich integral is developed in \S2 following the work of Le-Murakami
\cite\LMI\ and Kricker-Spence \cite\KS.  This leads to a lower bound for the power
series order of the Ohtsuki polynomial of a link, which is the essential
topological input to the bounds obtained in this paper.  In \S3 we study the
cyclotomic images of the Laurent polynomials introduced in \S1, and lower bounds
for their $p$-orders are given which strengthen previous bounds of Murakami
\cite\MI\ and Ohtsuki \cite\OI.  In \S4, these results are combined with a
diagonalizing lemma of Murakami and Ohtsuki \cite\MII\ to establish bounds on
the quantum $p$-orders of $3$-manifolds.  Finally in \S5 the examples needed to
sharpen these bounds are presented.

\heading{\S1. Laurent polynomials} \endheading

In this section we introduce two families of Laurent polynomials that arise in the
study of quantum invariants, one coming from number theory -- {\it quantum
integers} -- and the other from topology -- {\it Jones polynomials}.  These are
combined at the end of the section to form the {\it $p$-bracket} of a framed link,
which is the key ingredient in the surgery definition of the quantum invariant
$\tau_p$ given in \S4.

\section{The ring $\L$ and order}

Throughout the paper, $\L$ will denote the ring $\Z[t,t^{-1}]$ of integer Laurent
polynomials in an indeterminant $t$.  The variables $s=t^2$ and $q=t^4$ will also 
be used, as is common in the quantum topology literature.  These variables will
reappear in \S3 as roots of unity.

The following notion, and variations thereof, will be central to our
investigations.

\defn{1.1}
The {\it order} $\o(f)$ of a Laurent polynomial $f \in \L = \Z[t,t^{-1}]$ is the
order of $t=1$ as a zero of $f$.
\enddefn

\noindent
Thus $\o(f)$ is the lowest degree appearing in the Taylor expansion of $f$
about 1.  The related notion of {\it $p$-order}, for any prime $p$, is introduced
in \S3.  Lower bounds for the orders and $p$-orders of some of the
polynomials introduced in this section will be derived in later sections.

\remark{1.2}  
An equivalent way to define $\o(f)$ is by means of the substitution $h=t-1$.  This
embeds $\L$ into the ring of formal integer power series in $h$ (by mapping $t$ to
$1+h$ and $t^{-1}$ to $(1+h)^{-1} = 1-h+h^2-+\cdots$) and then $\o(f)$ is the smallest
power of $h$ with a non-zero coefficient in the series for $f$.  

It should be noted that the variable $h$ has often been used in the literature to
stand for $q-1$ or $\log(q)$ rather than $t-1$.  Thus Laurent polynomials in $t$
are transformed by substituting  $t=(1+h)^{1/4}$ or $t=\exp(h/4)$ into {\sl
rational} power series in $h$.  Fortunately the induced order valuations on $\L$ are
equal.  Indeed the lowest order terms for any substitution of the form $t =
1+ah+O(h^2)$, with $a\ne0$, will occur in the same degree, and so the order of $f$ is
well defined independent of such choices for the variable $h$.  To avoid confusion we
shall stick with the assignment $h=t-1$ throughout this paper, using different symbols
for other substitutions as the need arises.  For example $\hbar = \log(q)$ is used in
\S2.  Conway's substitution $z = s-s^{-1}$ is also useful in simplifying the formulas
in \S3, and provides the same notion of order since $z = 4h+O(h^2)$.

\section{Quantum integers}

For each integer $k$ define the {\it quantum integer} $[k] = (s^k-s^{-k}) /
(s-s^{-1})$, and more generally the {\it framed quantum integers}
$$
\ak = t^{a(k^2-1)} \, [k] 
\tag 1
$$
for any integer $a$; note that $(0,k]=[k]$.  These elements of $\L$ are ubiquitous
in the theory of quantum invariants.  Other versions of the quantum integers also
arise frequently, such as the classical {\it Gauss polynomials} $\br k =
(t^k-1)/(t-1)$, but we shall make little use of them (except briefly in \S3).  

Observe that $[k]$ can be written as a polynomial in $[2]$. Indeed it follows by
induction from the elementary identity $[k] = [2][k-1]-[k-2]$ that
$$
[k] = \sum_{j=0}^{k/2} (-1)^j {k-j-1\choose j} \, [2]^{k-2j-1}
\tag 2
$$
for $k\ge 0$ (these are renormalized Chebyshev polynomials), and $[-k]=-[k]$. 
The sum is (by convention) over all {\sl integers} \, $0\le j\le k/2$,
and so the upper limit is actually $k/2-1$ for $k$ even, since the binomial
coefficient vanishes when $j=k/2$, and $(k-1)/2$ for $k$ odd. 

More generally consider the two parameter family of {\it cabled quantum integers}
$$
\kc = {\sum_{j=0}^{k/2} (-1)^j {k-j-1\choose j} {k-2j-1\choose c} \,
[2]^{k-2j-1-c}}
\tag 3
$$
in $\L$; note that $[k,0)=[k]$.  These are defined for non-negative integers $k$
and $c$ by (3), and for arbitrary $k$ and $c$ by declaring
$[-k,c) =-[k,c)$ and $[k,c) = 0$ for $c<0$. They arise in conjuction with the
framed quantum integers in the formulas below for quantum invariants.  In
particular, the {\it $p$-sums} 
$$
(a|c) = \sum_{k=1}^{p/2} \ak\kc, 
\tag 4
$$
defined for any odd integer $p\ge 3$, play a special role.  Note that the only
dependence on $p$ is in the upper limit of summation, which is effectively
$(p-1)/2$ since $p$ is odd.  Since the value of $p$ is generally fixed, this
dependence is not made explicit in the notation.  

It will be shown in Proposition~3.7 that if $p$ is prime and $a$ is a multiple of
$p$, then the $p$-sum $(a|c)$ is divisible by $h^{p-3-2c}$ when viewed as an
element of the cyclotomic quotient $\lp$ of $\L$.  This is the main new technical
result used here to establish bounds for the quantum $p$-orders of $3$-manifolds.

\section{Jones polynomials}

The key topological input in the construction of quantum invariants of
$3$-manifolds is the {\it Jones polynomial} $V_L$ of an oriented link $L$
in the $3$-sphere \cite\Jo.  The version used here is a Laurent polynomial in the
variable $s=t^2$ characterized by $s^2V_{L_+} - s^{-2}V_{L_-} = (s-s^{-1})V_{L_0}$
and $V_{\bigcirc} = [2]$.   Here $L_+,L_-,L_0$ is the usual skein triple and
$\bigcirc$ is the unknot.  By convention $V_{\emptyset}=1$, where $\emptyset$ is
the empty link.

We shall actually use a variant $J_L$ of $V_L$, arising naturally
in the quantum group approach to the subject \cite\RTI \, \cite\KMI, which is
independent of the orientation on $L$. It is defined by 
$$ 
J_L = s^{3\lambda_L} V_L
\tag 5
$$ 
where $\lambda_L$ denotes the sum of all the pairwise linking numbers of $L$. 
Equivalently $J_L = (-1)^\ell K_L$, where $\ell$ is the number of components in
$L$ and $K_L$ is the Kauffman bracket of any zero-framed diagram for
$L$, normalized to be $1$ on the empty link \cite\Ka.  Clearly
$J_{\bigcirc} = [2]$, and more generally $J_{\bigcirc^c} = [2]^c$ where
$\bigcirc^c$ denotes the unlink of $c$ components, since $J$ is multiplicative
under distant unions $\sqcup$ of links.  (This polynomial is the zero-framed
version of the invariant $J_L$ defined in \S4 of \cite\KMI; the framed version has
value $(a,2]$ on the $a$-framed unknot.)

Using a cabling operation one may extend the definition of $J$ to the class of
{\sl colored} links $(L,k)$ in $S^3$. The {\it coloring} $k=(k_1,\dots,k_\ell)$
is a list of positive integers assigned to the components of $L$ (representing the
dimensions of simple modules in the quantum group approach). The {\it colored
Jones polynomial} of $(L,k)$ is then defined by the following formula, written in
the multi-index notation of [KM1] 
$$ 
J_{L,k} = \sum_{j=0}^{k/2} (-1)^j {k-j-1\choose j}
J_{L^{k-2j-1}}, 
\tag 6
$$ 
(This is the zero-framed version of the invariant $J_{L,\bold k}$ in \cite\KMI.)  
Here $L^c$ denotes the {\sl $c$-cable} of $L$, where the {\it cabling}
$c = (c_1\cdots,c_\ell)$ is a list of integers, obtained by replacing the $i$th
component of $L$ by $c_i$ parallel copies of itself with pairwise linking numbers
equal to zero; by convention $J_{L^c}=0$ if any $c_i<0$.  The multi-index notation
in (6) is to be interpreted as follows: the sum is over all lists
$j=(j_1,\cdots,j_\ell)$ with $0\le j_i<k_i/2$, and the signs and binomial
coefficients are products of the corresponding terms for each $i$.  Note that
$J_{L,k}$ can be defined for arbitrary integer colorings by requiring it to be an
odd function of any given color.

The reader may have noticed a similarity between equations
(2) and (6). Indeed, setting $k=2$ in (6) gives $J_{L,2} = J_L$, and so
(6) says that $J_{L,k}$ is the $k$th Chebyshev polynomial in $J_{L,2}$ with
cables replacing powers. In particular, taking the unknot for $L$ shows that
$J_{\bigcirc,k} = [k]$.  (The framed version of this invariant has
value $\ak$ on the $k$ colored unknot with framing $a$.)

\section{Ohtsuki polynomials}

Finally, we introduce another version $\phi_L$ of the Jones polynomial due
to Tomotada Ohtsuki \cite\OI.  First consider the free $\Z$-module
$\Cal L$ with basis consisting of all oriented links in $S^3$ (the orientation
is only relevant to the discussion of the Kontsevich integral below).   Note
that the invariant $J$ can be extended uniquely to a linear functional 
$J : \Cal L \to \L.$ 

Next consider the projection $\pi:\Cal L \to \Cal L$ defined by
$$
\pi(L) = \sum_{S<L}(-1)^{\ell-s} \, S|L
\tag 7
$$
where $S|L = S \sqcup \bigcirc^{\ell-s}$, obtained from $L$ by replacing each
component not in $S$ with a distant unknot.  The sum is over all sublinks $S$ of
$L$, including the empty link, and $\ell = \#L$, $s = \#S$.  That $\pi$ is a
projection follows readily from the fact that $\pi(S\sqcup T) = \pi(S)
\sqcup \pi(T)$ and $\pi(\bigcirc)=0$.  

Now define $\phi$ to be the composition $J\pi$,
that is
$$
\phi_L = J_{\pi(L)} = \sum_{S<L}[-2]^{\ell-s}J_S.
\tag 8
$$
Evidently $\phi_L$ is an integer Laurent polynomial in $s=t^2$, and in
particular an element of $\L$.  It will be called the {\it Ohtsuki polynomial}
of $L$ (see Remark 1.3 below).  

Observe that $J$ can be expressed in terms of $\phi$ as follows
$$ 
J_L = \sum_{S<L} [2]^{\ell-s} \, \phi_S. 
\tag 9
$$
To see this, note that any link $L$ can be recovered from the projections of its
sublinks as $L = \sum_{S<L} \pi(S) \sqcup \bigcirc^{\ell-s}$.  Indeed the right hand
side is $\sum_{R<S<L}(-1)^{s-r} \, R|L$ by definition, and this can be rewritten as
$\sum_{R<L} \,(\sum_{k=0}^{\ell-r} (-1)^k {\ell-r \choose k}) \, R|L$, which equals
$L$ since the inner sum vanishes for $R\ne L$.  Equation (9) follows.

More generally, the colored Jones polynomial $J_{L,k}$ can be rewritten as a
linear combination of $\phi_{L^c}$ for cables $c<k$ (i.e.\ $c_i<k_i$ for each
$i$), with coefficients the multi-index versions $[k,c) = \prod[k_i,c_i)$ of the
cabled quantum integers defined in (3): 

\lemma{1.3}
$J_{L,k} = \dsize\sum_{c=0}^{k-1} [k,c) \, \phi_{L^c}$
\endlemma

\proof  First use (9) to replace $J_{L^{k-2j-1}}$ in the definition (6)
of $J_{L,k}$ with the sum of $[2]^{k-2j-1-s}\phi_S$ over all sublinks $S$ of the
cabling $L^{k-2j-1}$ of $L$. Now each $S$ is again a cabling $L^c$ of $L$,
appearing ${k-2j-1 \choose c}$ times in the sum. Collecting terms yields the
result.
\endproof

\remark{1.4}
The polynomial $\phi_L$ defined here differs from Ohtsuki's original
version $\Phi_L$ by a factor of $[2]^{\ell}$.  Indeed by definition
$\Phi_L = (-1)^\ell X_{\delta(L)}$, where $X_S = J_S/[2]^s$ is Hitoshi Murakami's
normalization of $J_S$ [M1] and $\delta$ is the involution on $\Cal L$ given by 
$$
\delta(L) = \sum_{S<L}(-1)^s S.
$$
(That $\delta$ is an involution follows by an argument analogous to the
derivation of (9) above.)  Thus $\Phi_L = \phi_L/[2]^{\ell}$.  One advantage of
$\phi$ over $\Phi$ is that it is an honest Laurent polynomial, thus justifying its
name, taking values in $\L$ rather that $\L$ localized at $[2]$.  Ohtsuki's
normalization on the other hand leads to a pleasing symmetry between $X$ and
$\Phi$, namely $X_L = (-1)^\ell \Phi_{\delta(L)}$.  This is immediate from the
fact that $\delta$ is an involution, and yields an alternative derivation of
(9).

\section{The $p$-bracket}

Fix an odd integer $p\ge3$.  The key ingredient in the surgery definition of the
quantum invariant $\tau_p$ is the {\it $p$-bracket} of a {\sl framed} link
$L$ in the $3$-sphere.  It is an integer Laurent polynomial in $t$, that is an
element of $\L$, defined by 
$$
\br L = \sum_{k=1}^{p/2} (a,k]\, J_{L,k} 
\tag 10
$$ 
where $a = (a_1,\dots,a_\ell)$ is the list of integer framings on
$L$.  (As with the $p$-sums defined in (4), the dependence on $p$ is not explicit
in the notation, but is to be understood from the context.)  Multi-index notation
is being used as usual: the sum is over all colorings $k=(k_1,\dots,k_\ell)$ of $L$
with $1 \le k_i < p/2$ for each $i$, and the coefficients $(a,k]$ are multi-index
versions $\prod (a_i,k_i]$ of the framed quantum integers defined in (1).  

For any integer $a$, set $b_a = \br{\bigcirc_a}$, the $p$-bracket of the $a$-framed
unknot.  This coincides with the associated zero-cabled $p$-sum defined in (4)
$$
b_a = (a|0)
\tag 11
$$  
since $J_{\bigcirc,k}=[k]$.

We conclude this section with an expression for the $p$-bracket of $L$ in terms of
Ohtsuki polynomials of cablings $L^c$ of $L$ and the associated multi-index
versions $\ac = \prod (a_i|c_i)$ of the $p$-sums defined in (4).

\proposition{1.5}  The $p$-bracket of a framed link $L$ in the
$3$-sphere with framings $a$ can be written as
$
\br L = \sum_{c=0}^n \ac \, \phi_{L^c}
$
where $n=(p-3)/2$.
\endproposition

\proof By Lemma 1.3 and the definition of the $p$-bracket,
$$
\br L = \sum_{k=1}^{p/2} \,  \sum_{c=0}^{k-1} \ak\kc \, \phi_{L^c}.
$$
Since $\kc = 0$ for $c\ge k$, the upper limit of the inner sum can
be replaced with $n$.  Switching the order of summation then gives the result.
\endproof

%%%%%%% section 2 %%%%%%%

\heading{\S2. Relations with the Kontsevich integral} \endheading

The Jones and Ohtsuki polynomials of a link in the $3$-sphere can both be
interpreted in terms of the Konsevich integral of the link.  For the Jones
polynomial, this was first made explicit in the work of Le and Murakami
\cite\LMI, motivated by earlier work of Drinfeld on quasi-Hopf algebras \cite\Dr. 
For the Ohtsuki polynomial, this was elucidated in a recent paper of Kricker and
Spence \cite\KS, following ideas suggested in the seminal work of Le on
finite type invariants of homology spheres \cite\Le.  This interpretation of the
Ohtsuki polynomial can be made particularly transparent using the projection
$\pi:\Cal L \to \Cal L$ defined in (7) , as explained below.

\section{The Kontsevich integral}

Recall the {\it Kontsevich integral} $\hat Z_L$ of a {\sl framed} oriented
link $L$ in $S^3$, normalized as in \cite\LMO.  It is an element of the
completion $\A(L)$ (with respect to degree) of the rational vector space generated
by Feynman diagrams on $L$ modulo the appropriate relations (AS, IHX and STU). 
Here a {\it Feynman diagram} on $L$ consists of an abstract vertex-oriented
uni-trivalent graph with all of its univalent vertices on $L$ (see \cite\LMO\ for
details).  This graph is generally referred to as the {\it dashed graph} of the
diagram, as it is drawn with dashed lines in pictures.  The univalent vertices are
called {\it external} vertices, and the trivalent ones are called {\it internal}. 
The {\it degree} of the diagram is half the total number of vertices.   

Any element in $\A(L)$ can be expanded (in many ways) as an {\sl infinite} linear
combination of diagrams; such an expansion will be called a {\it Feynman series}
for the element.  A simply connected component of the dashed graph of any diagram in
the series will be called a {\it tree}, and if it has all but possibly one external
vertex lying on a single component of $L$, then it will be called a {\it thin tree}. 
These notions will be used later in this section.

For the present purposes we shall only consider {\sl zero-framed} links, and shall
also assume that all links come equipped with an ordering for their components so
that $\A(L)$ is identified with $\A(\bigcirc^\ell)$ for any $\ell$-component link
$L$.  The Kontsevich integral can then be viewed as a $\Z$-linear map 
$
\hat Z : \Cal L \to \A
$
where where $\Cal L$ is the free $\Z$-module generated by oriented links in
$S^3$ and $\A = \oplus_{\ell=0}^{\infty} \A(\bigcirc^\ell)$.  This map is a
morphism for distant unions, i.e.\ $\hat Z_{J\sqcup K} = \hat Z_J \sqcup \hat
Z_K$.  As is standard practice, $\nu$ will denote the value of $\hat Z$ on the
unknot, which is a unit in $\A(\bigcirc)$ with respect to the connected sum
operation $\#$ (see
\S1.1 in \cite\LMO).  

In fact it is convenient to use a related invariant $Z : \Cal L \to \A$ defined by
$$
Z_L = \hat Z_L \, \# \, (\nu^{-1})^\ell
\tag 12
$$
meaning that one takes the connected sum of a copy of $\nu^{-1}$ with each
component of $L$ in each term of $\hat Z_L$.  (This is the zero-framed version of
the invariant $\hat Z'$ of \cite\KS.)  Note that $Z$ is also a morphism for
distant unions.

\section{Connections with the Jones and Ohtsuki polynomials}

To interpret the Jones and Ohtsuki polynomials (which take values in the
ring $\L = \Z[t,t^{-1}]$) in terms of the Kontsevich integral, it is useful to
expand these polynomials as power series in $\hbar = \log(q)$ where
$q=t^4$ (see Remark 1.2).  In particular, consider the embedding $E:\L \to
\Q[[\hbar]]$ which maps $t$ to $\exp(\hbar/4)$.  For any function $F$ mapping into
$\L$, write $\hat F$ for the composition $EF$.  Thus for example $\hat J_L$ is the
power series obtained from the Jones polynomial $J_L$ by substituting
$\exp(\hbar)$ for $q$, and $\hat\phi_L = \hat J_{\pi(L)}$ is the series for the
corresponding Ohtsuki polynomial.

Now following \cite\KS, let $\G : \A \to \Q[[\hbar]]$ denote the weight system
associated with the trace form in the fundamental representation of $sl(2,\C)$. 
In fact, it is convenient to use a renormalization $W$ of $\G$ defined by
$
W(D) = (\widehat{[2]}/2)^\ell \, \G(D)
$
on any Feynman diagram $D$ on an $\ell$ component link.  Note that $W$ preserves
order (the lowest degree in a power series) since $\G$ preserves degree and
$\widehat{[2]} = 2 + \hbar^2/4 + \cdots$ is of order zero.  Also let $P:\A\to\A$
denote the projection which sends a diagram $D$ on $L$ to itself if it has
vertices on {\sl every} component of $L$, and to zero otherwise.  Then the
interpretation of the Jones and Ohtsuki polynomials in terms of the Kontsevich
integral can be expressed as follows.

\lemma{2.1}
The diagram
$$
\CD
\Cal L   @>\pi>>   \Cal L   @>J>>   \L   \\
@VZVV              @VZVV           @VVEV   \\
\A       @>>P>      \A      @>>W>   \Q[[\hbar]] 
\endCD
$$
commutes, and so \ {\rm (a) (Le-Murakami)} $\hat J = WZ$, and 
\ {\rm (b) (Kricker-Spence)} $\hat\phi = WPZ$.  
\endlemma

\proof 
The commutativity of the square on the right, which is equivalent to (a), is the
result of Le-Murakami \cite\LMI\ mentioned above.  Actually they show $EJ =
\G\hat Z$, but it is easy to verify that $\G\hat Z = W Z$ using the
identities $\G(D\# D')={1\over2}\G(D)\G(D')$, $\G(\bigcirc) = 2$ and
$\G(\nu)=\widehat{[2]}$ (see \cite\LMII\cite\KS).

For the square on the left, observe that for any sublink $S$ of a link $L\in\Cal L$,
the Kontsevich integral $Z_{S|L}$ of $S|L = S\sqcup\bigcirc^{\ell-s}$, which appears
in the definition of $\pi$ in (7), is the sum of the terms in $Z_L$ whose diagrams
have all their vertices on $S$.  This is an elementary consequence of the work of Le
and Murakami which describes how $\hat Z_S$ and $\hat Z_L$ are related (see
Proposition 1.1 in \cite\LMO).  Now the commutativity of the left-hand square
follows by the inclusion-exclusion principle of combinatorics.  This implies the
commutativity of the outermost rectangle (shown directly in \cite\KS) which is
equivalent to (b).
\endproof

\vskip -5 pt
\section{Orders of Ohtsuki polynomials}

Kricker and Spence's formula $\hat \phi_L = W(P(Z_L))$ led them to a striking lower
bound for the order of the Ohtsuki polynomial of any link $L$ whose pairwise
linking numbers all vanish.  Such a link will be called a {\it diagonal}
(or {\it algebraically split}) link.

\theorem{2.2} {\rm (Kricker-Spence \cite{\KS,\S3})}  Let $L$ be a diagonal link
with $\ell$ components.  Then the order $\o(\phi_L)$ of the Ohtsuki polynomial of
$L$ is greater than or equal to $4\ell/3$.
\endtheorem

Using different methods, Ohtsuki had previously obtained a bound in terms of the 
{\it maximum cabling index} of the link, which is by definition the maximum
number $m$ of mutually parallel, algebraically unlinked components of the link. 
Although weaker than the Kricker-Spence bound in many cases, Ohtsuki's bound is
stronger whenever $m>\ell/3$, a situation that will arise in \S4 in establishing
the $p$-order bounds for 3-manifolds of small betti number.

\theorem{2.3} {\rm (Ohtsuki \cite{\OI,\S3.4})}  Let $L$ be a diagonal link with
$\ell$ components and maximum cabling index $m$.  Then the order $\o(\phi_L)$ of the
Ohtsuki polynomial of $L$ is greater than or equal to $\ell+m$.
\endtheorem

The theorem of Kricker and Spence follows from their formula by an easy counting
argument, after making some elementary observations about the Kontsevich integral of
a diagonal link.  In fact a recent result of Habegger and Masbaum \cite\HM\
(motivated by earlier work of Rozansky \cite\RI\ and Le \cite\Le) shows how to
generalize these observations to links with vanishing higher order linking numbers,
and so we present the proof in this more general context.  Since Ohtsuki's theorem
can also be viewed in this context, we formulate a theorem below that includes both
results. 
  
First define the {\it Milnor degree} of a link $L$ in $S^3$ to be the degree
(= length$-1$) of the first nonvanishing $\bar\mu$-invariant of $L$ (see
\cite\Mi).  Thus {\sl every} link has degree $\ge1$, while the diagonal links are
those of degree $\ge2$.  If all the $\bar\mu$-invariants of $L$ vanish (e.g.\ if
$L$ is a knot, a boundary link or more generally the fusion of a boundary link
\cite\Co) then $L$ is said to have infinite Milnor degree.

Now the result of Habegger and Masbaum  (\S6.10 in \cite\HM) is that the Kontsevich
integral $\hat Z_L$ of a link $L$ of Milnor degree $d$ has a Feynman series in which
every tree has degree at least $d$.  Of course this is also true for the
normalization $Z_L$, since $\nu = \hat Z_{\bigcirc}$ can be written as a linear
combination of treeless diagrams (by the same result).  In fact the proof in
\cite\HM\ shows more:  One can arrange that for each sublink $S$ of $L$, every tree
with all of its external vertices on $S$ has degree no less than the Milnor degree
of $S$; thus for example there will be no chords between any pair of components
with linking number zero, even if $L$ is not diagonal.  Furthermore, thin trees
(ones with all but possibly one vertex on a single component of $L$) can be avoided
in the series if $L$ is diagonal.  A series satisfying these conditions will be
called a {\it Milnor series} for $Z_L$.

\lemma{2.4} The Kontsevich integral $Z_L$ of any link $L$ has a Milnor series.
\endtheorem

\proof
We assume that the reader is familiar with the arguments presented in \S6 of
\cite\HM, and adopt the terminology and notation used there.

Choose a string link $\lambda$ representing $L$.  By the proof of \S6.9--10 in
\cite\HM, it suffices to show that $Z_\lambda^t$ (which is an element of
$\A^t(\ell)$) has a Milnor series.  Now the key observation, as in \cite\HM, is
that $Z_\lambda^t$ is {\sl group like}, and so can be expressed as the exponential
of some primitive $\xi$ which is a linear combination of tree diagrams.  For each
string sublink $\sigma<\lambda$, let $\xi_\sigma$ denote the sum of the diagrams in
$\xi$ which have vertices on every component of $\sigma$, but on no other
components of $\lambda$.  Thus $\xi=\sum_{\sigma<\lambda}\xi_\sigma$.  

As in the proof of Lemma 2.1 above, it follows from a formula of Le and Murakami
(Proposition 1.1.3 in \cite\LMO) that $Z^t_{\sigma|\lambda} = \exp(\sum_{\tau<\sigma}
\xi_\tau)$, and so by Theorem 6.2 in \cite\HM, no tree in $\xi_\sigma$ is of degree
less that the Milnor degree of $\sigma$.  In particular if $L$ is diagonal, then
there are no chords, and any other thin trees can be omitted as well since they
vanish in $\A^t(\ell)$.  Therefore $\exp(\xi)$ is a Milnor series for $Z^t_\lambda$.
\endproof

\theorem{2.5}  Let $L$ be an $\ell$-component link of Milnor degree $d$.  Then the
order $\o(\phi_L)$ of the Ohtsuki polynomial of $L$ is greater than or equal to
$2\ell d/(d+1)$.  Furthermore, if $L$ is diagonal {\rm($d\ge2$)} with maximum
cabling index $m$, then $\o(\phi_L) \ge \ell+m$.    
\endtheorem

In particular every link satisfies $\o(\phi_L)\ge\ell$.  For diagonal links the
bounds can be incorporated into a single formula
$$
\o(\phi_L) \ge \ell + \max({d-1\over d+1}\ell,m),
$$
and for links of infinite Milnor degree $\o(\phi_L)\ge2\ell$. 

\pf{2.5} By Lemma 2.1 and the fact that $W$ preserves order, it suffices to show
that $P(Z_L)$ has a Feynman series consisting of diagrams of degree $\ge 2\ell
d/(d+1)$, and $\ge \ell+m$ when $L$ is diagonal.  

First observe that the degree $d_D$ of a diagram $D$ can be computed as the
difference $x_D-e_D$, where $x_D$ is the number of external vertices in $D$ and
$e_D$ is the euler characteristic of its dashed graph, or alternatively as the sum
of the {\it weights} of the external vertices of $D$.  Here the weight of a vertex is
defined as the ratio $d_C/x_C$, where $C$ is the {\sl component} of $D$ containing
the vertex.  By the first formula for the degree, this weight is $d/(d+1)$ if $C$ is
a tree of degree $d$, and is $\ge1$ otherwise.  

It follows from Lemma 2.4 that $Z_L$ has a Feynman series with all external
vertices of weight $\ge d/(d+1)$.  The corresponding series for $P(Z_L)$ is
obtained by eliminating certain diagrams, leaving only those with at least one
vertex on each component of $L$.  Moreover, diagrams with {\sl exactly} one vertex
on some component vanish in $\A$ by the STU relation.  The remaining diagrams have
at least $2\ell$ external vertices, and therefore degree $\ge 2\ell d/(d+1)$.

Now consider the case when $L$ is diagonal.  By hypothesis $L = J\cup K^m$
where $J$ has $j = \ell-m$ components and $K$ is a knot.  Applying Lemma 2.4 again,
choose a series for $P(Z_{J\cup K})$ consisting of diagrams with no thin trees, no
trees of degree $<d$, and with at least two vertices on each component of $J\cup
K$.  Cabling this series (following \cite{\LMII,\S4.1}) produces a similar series for
$P(Z_L)$ with the additional property that every tree has at least two vertices on
$J$.  Now any diagram $D$ in this series has degree $d_D = x_D-e_D \ge
x_D-t_D$, where $t_D$ denotes the number of trees in $D$.  Since $x_D = 2\ell +
k$ for some $k\ge0$,  there are at most $2j+k$ vertices on $J$, and so $t_D\le
j+k/2$.  Therefore $d_D \ge 2\ell+k-(j+k/2) \ge \ell+m$.  
\endproof

\heading{\S3. Cyclotomacy}\endheading

In this section we introduce {\sl cyclotomic orders} in the ring $\L =
\Z[t,t^{-1}]$ of Laurent polynomials.  In particular, for each prime $p$ we define
the {\it $p$-order} on $\L$ in terms of a certain distinguished valuation $\op$ on
the quotient ring $\lp$ by the ideal generated by the $p$th cyclotomic
polynomial.  Bounds will then be established for the $p$-orders of the polynomials
defined in \S1, namely the $p$-sums (Propositions 3.6 and 3.7), and the Ohtsuki
polynomials and $p$-brackets of links (Proposition 3.5 and Theorem 3.10). 

\section{The ring $\lp$ of cyclotomic integers}

Fix an odd prime $p=2n+3$.  The
cyclotomic polynomial $\varphi_p = (t^p-1)/(t-1)$ generates a prime ideal in $\L =
\Z[t,t^{-1}]$. Let $\lp$ denote the quotient of $\L$ by this ideal,
$$
\lp = \L/(\pc),
$$ 
and $\qp$ denote the field of fractions of $\lp$.  Of course $\qp$ can be
identified with the cyclotomic field of {\sl complex} $p$th roots of
unity, and $\lp$ with its ring of integers.  Indeed for any primitive complex
$p$th root $\zeta$, there is a unique ring homomorphism $\L \to \C$ mapping $t$ to
$\zeta$, and this induces an isomorphism $\qp\cong\Q(\zeta)$ carrying $\lp$ to
$\Z[\zeta]$.  Various properties of $\lp$ can be deduced from this observation. 
For example the classical fact that the ring of integers in a number field is
a Dedekind domain shows that ideals in $\lp$ factor uniquely
into prime ideals.

For convenience we retain the symbol $t$ for the image of $t$ in $\lp$, and
continue to use the notation $s=t^2$ and $q=t^4$.  Thus a Laurent polynomial in
$t$ may be used to represent either an element in $\L$ or an element in $\lp$,
depending upon the context.  Note however that these two elements may have very
different number theoretic properties.  For example $t$, $s$ and $q$ become
$p$th roots of unity, and thus units, in $\lp$.  The quantum integers $[k] =
(s^k-s^{-k}) / (s-s^{-1})$ are also units in $\lp$ if $k$ is prime to $p$, with
$[k]^{-1} = [k\bar k]/[k]$ (clearly an integral polynomial in $t$) where $\bar k$
is any mod $p$ inverse of $k$.  Similarly the Gauss polynomials $\br k =
(t^k-1)/(t-1)$ are units in $\lp$ for $k$ prime to $p$.  Note that $[k]$ is an odd
function of $k$ while $\br k$ is not.  Both are periodic of period $p$.

\section{The prime ideal $H$ and $p$-order in $\lp$}

The polynomial $h = t-1$, which is evidently a prime in $\L$, remains prime as
an element of $\lp$ (see for example \S1.4 in \cite\Wa).  In contrast $p$ is not a
prime in $\lp$.  It is in fact a power of $h$ times a unit, as seen by the
calculation $p = \varphi_p(1) = (t-1)\cdots (t^{p-1}-1) = \br{p-1}!\,h^{p-1}$. 
(Here $\br{p-1}!$ is the product $\br1\cdots\br{p-1}$ of Gauss polynomials.)  Thus
$h$ generates a prime ideal $H$ in $\lp$ and 
$$
P = H^{p-1}
$$
where $P=(p)$ is the ideal generated by $p$.  This is an instance of the prime
factorization of ideals in $\lp$.  The uniqueness of this factorization shows that
$H$ is the only prime ideal in $\lp$ containing $p$. 

The powers of $H$ form a descending sequence  
$$
\lp = H^0 \supset H^1 \supset H^2 \supset \cdots
$$
of ideals with trivial intersection.  Additively, each of these is free
abelian of rank $p-1$.  Indeed the consecutive powers $h^k,\cdots,h^{k+p-2}$ form a
basis for $H^k$ as a $\Z$-module.  Thus any non-zero element $x$ in $\lp$
lies in a unique smallest $H^{k_0}$, and can be written uniquely as an integer
polynomial $x_k = \sum_{d=0}^{p-2} x_{k,d}h^{k+d}$ in $h$ for each $k\le k_0$.  We
will call $x_{k_0}$ the {\it normal form} of $x$, and $x_0$ (which is of degree
$\le p-2$) the {\it reduced form} of $x$.   

The integer $k_0$, which is the {\sl order} of the normal form of $x$ in the sense
of Remark 1.2, is called the {\it $p$-order} of $x$, and will be denoted
by $\op(x)$.  In simple terms $\op(x)$ is the exponent of the highest power of $h$
that divides $x$ (as an element of $\lp$).  This notion can be lifted to the ring
$\L$ by declaring the $p$-order of a Laurent polynomial $f$, also written $\op(f)$,
to be the $p$-order of its image in $\lp$.  These notions will be recast
below in the context of valuation theory.

\remarks{3.1}
(a)  In discussing $p$-order, $h$ may be replaced by any other
generator of $H$, that is any associate of $h$ in $\lp$.  Among the possible
choices are $q-1$ and $z = s-s^{-1}$ (cf.\ Remark 1.2), or more generally 
any of the elements $h_{j,k} = t^j-t^k$  for $j\not\equiv k\pmod p$, since $h_{j,k}
= (\br{j\!-\!k}t^k) \, h$.  Furthermore some formulas are more revealing with a
different choice of generator.  For example the formula $p = \br{p-1}!h^{p-1}$
above has an analogue using the generator $z$, namely $p = [p-1]!z^{p-1}$.  This
is proved in the same way, but starting with the factorization $p =
(q-1)\cdots(q^{p-1}-1)$.  But now exploiting the fact that $[k]$ is an {\sl odd}
periodic function of $k$, this can be rewritten as $p = (-1)^m([m]!z^m)^2$ where
$m=(p-1)/2$.  This recovers the well known fact that $(-1)^mp$ is a square in
$\lp$.  Moreover, its square roots can be identified with the Gauss sums
$$
G_a \underset \,\text{def} \to = \,\sum_{k=1}^p t^{ak^2} =
\fracwithdelims()ap(-1)^m[m]!z^m
$$
by an easy argument using the classical formula $\prod_{k=1}^m
(t^{2k-1}-t^{-(2k-1)})$ for $G = G_1$ \cite{\IR,\S6.4.4}.  Here $(\!-\!)$ is the
Legendre symbol.

\smallskip

(b) An integral polynomial $\sum_{d=0}^{p-2}x_d h^{k+d}$ in $h$ with $x_0\ne0$ is in
normal form if and only if $x_0$ is prime to $p$.  For if $p$ divides $x_0$, then the
first term $x_0h^k$ can be killed by substracting a multiple of the cyclotomic
polynomial $\pc$, since $\pc$ has constant coefficient $p$ when written as a
polynomial in $h$; indeed $\pc = (t^p-1)/(t-1) = ((h+1)^p-1)/h$, and so
$$
\pc =  p + \sum_{k=2}^{p-1}
{p\choose k} h^{k-1} + h^{p-1}
\tag 13
$$
by the binomial theorem.  The converse follows from the fact that $p$ is the only
prime (rational) integer divisible by $h$.

\smallskip

(c) If an element of $\lp$ is given as the image of a Laurent polynomial
in $t$, then it can be put in either normal or reduced form in the following way: 
First rewrite the polynomial using the relation $t^p=1$ as an {\sl honest}
polynomial in $t$, i.e.\ an element in $\Z[t]$.  Then substitute $t=1+h$ to
obtain an element in $\Z[h]$.  Finally subtract a suitable multiple of
$\varphi_p$ to put this in the required form, working from the ``bottom-up"
for normal form as indicated in the previous remark, and from the ``top-down" for
reduced form.  For example if $p=3$, then working mod $\varphi_3 = 3+3h+h^2$ we
have $t+2t^{-1}
\equiv t+2t^2 \equiv 3+5h+2h^2 \equiv 2h+h^2$ (normal form) or $-3-2h$ (reduced
form).  Thus $\o_3(t+2t^{-1}) = \o_3(2h+h^2) = 1$.

\section{Valuation theory}

A {\it prevaluation} on an integral domain $D$ is a nonconstant map
$\v:D\to\Z\cup\{\infty\}$ satisfying
$$
\align
\text{(a)} \quad &\v(xy) = \v(x)+\v(y) \\
\text{(b)} \quad  &\v(x+y) \ge \min(\v(x),\v(y)).
\endalign
$$  
It follows from (a) that $\v(1)=0$, and so
$\v(u^{-1})=-\v(u)$ for any unit $u$ in $D$. In particular $\v(x/y)=\v(x)-\v(y)$
if $D$ is a field. It also follows from (a) that $\v(0)=\infty$. 

The set $\Cal R_\v$ of elements $x$ in $D$ with $\v(x)=\infty$ is called the {\it
radical} of $\v$, and is evidently an ideal in $D$. If $\Cal R_\v=0$, then $\v$
is called a {\it valuation} (and so any prevaluation on a field is {\sl a
fortiori} a valuation). The set $\calo_\v$ of elements with $\v(x)\ge0$ is a
subring of $D$ called the ({\it pre}){\it valuation ring} of $\v$.
If $\calo_\v=D$ then $\v$ is said to be {\it positive}.  More generally consider
the descending filtration $\cdots \supset \calodv \supset \calo_{\v}^{d+1} \supset
\cdots$ of subrings of $D$, where  $\calodv = \{x\,|\,\v(x) \ge d\}$.  In
particular $\calo^{\infty}_\v = \Cal R_\v$ and $\calo^0_\v=\calo_\v$. These are
all ideals when $\v$ is positive. In general, the value $\v(x)$ is determined by
the position of $x$ in this filtration; $\calodv - \calo^{d+1}_\v$ is
the set of elements with value exactly $d$.

The standard example of a valuation in elementary number theory is the
{\it $p$-adic valuation} $\vp$ on $\Q$ which assigns to each rational number $x$
the exponent of $p$ in the prime decomposition of $x$.  It restricts to a positive
valuation on the integers.  The $p$-adic valuation has a natural generalization to
the cyclotomic field $\qp$:

\section{The $p$-adic valuation and $p$-order in $Q_p$}

The {\it $p$-adic valuation} $\op$ on $\qp$ assigns to each element $x$ the
exponent of the prime ideal $H$, the unique ideal in $\lp$ lying over $p$, in the
prime decomposition of the fractional ideal generated by $x$.  (This is an
instance of the general construction of valuations on the field of fractions of a
Dedekind domain $D$ from prime ideals in $D$ \cite\BS.)   The integer $\op(x)$ will
be called the {\it $p$-order} of $x$.  This clearly coincides with the definition
given above for $x \in \lp$, since $\calo_{\op}^d = H^d$, and so in this case
$\op(x)$ can be computed as the order of the normal form of $x$ (which allows for
computation in general since $\op(x/y) = \op(x)-\op(y)$).  Alternatively, the
$p$-order of $x\in\lp$ can be computed from the reduced form $\sum_{d=0}^{p-2}
x_dh^d$ of $x$ by   
$
\op(x) = \min_d(\op(x_dh^d)) = \min_d((p-1)\vp(x_d) + d)
$
(the first equality holds since the terms in the sum have {\sl distinct}
$p$-orders).  

\remark{3.2}
The last computational scheme suggests introducing the $\Z$-linear projections
$$\pi^d:\lp\to\Z_{p^k}, \quad x \to x_{\bar d} \text{ (mod }p^k)
$$
for $d\ge0$.  Here $x_{\bar d}$ is the coefficient of $h^{\bar d}$ in the reduced
form of $x$ and $k = 1+\lfloor d/(p-1) \rfloor$, where $\bar d$ is the least
residue of $d$ mod $(p-1)$ and $\lfloor\ \rfloor$ is the greatest integer
function.  Now $\op(x)$ can be defined as the smallest integer $d$ for which
$\pi^d(x)\ne0$.  The motivation for this point of view comes from the study of the
``$p$-adic completion" $\hat\L_p = \varprojlim \lp/P^k$ of the ring $\lp$. 
Indeed, the $\pi^d$'s can be grouped in strings of length $p-1$ to give
projections $\pi_k:\lp\to\lp/P^k \cong (\Z_{p^k})^{p-1}$ which show
$\hat\L_p$ isomorphic to a direct sum of $p-1$ copies of the $p$-adic integers. 
Note that $\lp$ embeds in $\hat\L_p$, i.e.\ any element $x$ in $\lp$ can be
recovered from its projections $\pi^d(x)$ for $d=0,1,2,\dots$.  
\endremark

\vskip -5 pt
\section{$p$-order in $\L$}

The restriction of the $p$-adic valuation to $\lp$ induces a positive {\sl
pre}valuation on $\L$, 
$
\op : \L \to \Z\cup\{\infty\},
$
by composing with the natural projection $\L\to\lp$.  This will be called the {\it
$p$-adic prevaluation} on $\L$.  It has radical $(\pc)$, and more generally
$\calo^d_\op = (h^d,\pc)$ where $h^{\infty}=0$ by convention.  The value
$\op(f)$ will be referred to as the {\it $p$-order} of $f$.  In other words:

\defn{3.3}
The {\it $p$-order} $\op(f)$ of a Laurent polynomial $f \in \L$ is the
$p$-order of the image of $f$ in $\lp$.  
\enddefn

Note that $\op(f)$ depends only on the equivalence class of $f \pmod{\pc}$, and so
its computation is facilitated by an appropriate choice of representative. 
Typically one uses either the normal or reduced form of $f$, viewed as an element
of $\lp$ (see Remark 3.1 and the computational schemes for $p$-order in $\lp$
discussed above).  One can, however, glean some information about $\op(f)$ without
reducing $f$ (i.e.\ avoiding the last step of the process described in Remark
3.1c).  In particular, lower bounds on $\op(f)$ can be obtained by comparing $\op$
with two other naturally defined (pre)valuations on $\L$, which we discuss next.

\section{Order and mod $p$\,-\,order in $\L$}

The {\it order valuation} $\o$ on the power series ring $\Z[[h]]$, defined by 
$$
\o(\tsize\sum x_d h^d) = \min\{d\,|\,x_d\ne0\},
$$  
induces a positive valuation on $\L$ by means of the usual embedding $t\mapsto
1+h$.  We call $\o(f)$ the {\it order} of $f$ (cf.\ Definition 1.1).  Similarly
define the prevaluation $\osp$ on $\Z[[h]]$ by 
$$
\osp(\tsize\sum x_d h^d) = \min\{d\,|\,x_d\not\equiv0\ \text{mod $p$}\}.
$$  
This induces a positive prevaluation on $\L$ with radical $(p)$, and $\osp(f)$ is
called the {\it mod $p$\,-\,order} of $f$.  The three prevaluations $\o$, $\op$
and $\osp$ on $\L$ are related as follows (cf.\ \cite{\MI,\S5.5}\ and
\cite{\OI,\S7.3}).

\lemma{3.4}
Let $f$ be a Laurent polynomial in $t$ and $f'$ denote its derivative with
respect to $t$ (or equivalently with respect to $h=t-1$).  Then for any
integer $d<\o(f)+p$,
\smallskip 
{\rm (a)} $\op(f)\ge\o(f)$, $\osp(f)\ge\o(f)$
\smallskip
{\rm (b)} $\op(f)\ge d \iff \osp(f)\ge d$
\smallskip
{\rm (c)} $\op(f)\ge d \implies \op(f')\ge d-1$.
\smallskip
\endlemma

\smallskip
Note that (a) and (b) say that $\op$ and $\osp$ are both $\ge\o$, and that
$\op = \osp$ whenever either one is $\le\o+p-1$.  In fact there are no further
relations among these prevaluations.  This can be seen using the examples
$h^i+h^jp+h^k\pc$ with $i\le j,k$, which have $\o=i$, $\op=j+p-1$ and
$\osp=k+p-1$.

\pf{3.4}
First observe that $f$ is a power of $t$ times a polynomial in $\Z[t]$ with
nonzero constant term.  Since powers of $t$ are units in $\L$, which have value $0$
with respect to any positive prevaluation on $\L$, it suffices to prove the lemma
for $f\in\Z[t]$ with $\o(f)=0$.  In the language of ideals, it says \  (a) $\calo^d
\subset \calodp$ for all $d$, \ (b) $\calodp = \calodsp$ for $d<p$, and \  (c)
$(\calodp)' \subset \calo^{d-1}_p$ for $d<p$, where for brevity we denote the
ideals $\calo_{\o_{\!-}}^d$ by $\calo_{\!-}^d$.

Since $\Z[t]$ is identified with the polynomial ring $\Z[h]$ under the substitution
$t=1+h$, the lemma is really a statement about the prevaluations $\o$, $\op$ and
$\osp$ on $\Z[h]$, where $\calo^d = (h^d)$, $\calodp = (h^d,\pc)$ and $\calodsp =
(h^d,p)$.  Now (a) is obvious and (c) follows from (b) using the elementary
observation that the derivative of any polynomial in $\calodsp$ is in
$\calo^{d-1}_{\!/\!p}$, since $d<p$.  For (b) it suffices to show
$\pc\in\calo^{p-1}_{\!/\!p}$ and $p\in\calo^{p-1}_{\!p}$.  But this follows
immediately from the observation that $\pc$ and $p$ are associates in the ring
$\Z[h]/(h^{p-1})$, which can be seen as follows:  By (13) $\pc = pu+h^{p-1},$
where $u = {1\over p}\sum_{k=1}^{p-1} {p\choose k} h^{k-1},$ and evidently $u = 1
+ O(h)$ is a unit in
$\Z[h]/(h^{p-1})$.
\endproof

\vskip -5 pt
\section{Bounds for the $p$-order}

As a trivial consequence of Theorem 2.5 and Lemma 3.4a we obtain lower bounds for
the $p$-orders of Ohtsuki polynomials of links in the $3$-sphere.

\proposition{3.5} Let $L$ be an $\ell$-component link of Milnor degree $d$.  Then
for any odd prime $p$,
$$
\op(\phi_L) \ge 2\ell d/(d+1).
$$
Furthermore, if $L$ is diagonal with maximum cabling index $m$, then
$\op(\phi_L)\ge\ell+m$.
\endtheorem

The following result of Ohtsuki \cite{\OI,\S7.2} gives bounds on the $p$-orders of
the $p$-sums $\ac$ defined in (4), for arbitrary integers $a$ and $c$ (see also
\cite{\MI,\S5.4}).

\proposition{3.6} {\rm (Ohtsuki)}  Let $p=2n+3$ be an odd prime. 
Then \ $\op \ac \ge n-c$.
\endproposition

This result can be strengthened when $a \equiv 0 \pmod p$ using Lemma 3.4.  

\proposition{3.7}
Let $p=2n+3$ be an odd prime and $a$ be a multiple of $p$.  Then \ $\op \ac \ge
2(n-c)$.
\endproposition

\proof
First observe that $\ac$ depends only on the residue class of $a\!\pmod p$, and so
$\ac = (0|c)$. For convenience write $(c)$ for $2(a|c)$ and set $m=n+1=(p-1)/2$.
Thus 
$ 
(c) = \sum_{k=-m}^m [k]\kc,
$
and it suffices to prove $\op(c) \ge 2(n-c)$ since $\op(2)=0$. If $c<0$, then
$\kc=0$ and so $\op(c)=\infty$. If $c>n$, then $2(n-c)<0$ and there is nothing
to prove, since $\op$ is positive. Thus we assume $0 \le c \le
n$.

Now consider the more general sum
$
(j,c) = \sum_{k=-m}^m [k]^{(j)}\kc
$
where $^{(j)}$ denotes the $j$th derivative (with respect to $t$ or
$h$). We will show
$$
\op(j,c) \ge 2(n-c)-j
\tag 14
$$
(for $0 \le c \le n$) by double induction on $c$ and $j$.  The proposition
follows by setting $j=0$. 

\subsection{Initial step of the induction} 

It must be shown that $\op(j,0)\ge 2n-j$, where $(j,0) = \sum_{k=-m}^m
[k]^{(j)}[k].$  By Lemma 3.4b it suffices to prove
$$
\osp(j,0)\ge 2n-j.
\tag 15
$$

To facilitate the proof, consider the ring $P=\Q[\cdot]$ of polynomials with
rational coefficients in an un-named variable $\cdot$ and let $\po$ be the subring
of all odd polynomials.  Set $\P = P[[h]]$ and $\ppo = \po[[h]]$. For
each integer $r\ge0$, let $P^r$ denote the subring of $P$ consisting of all
polynomials of degree $\le r$ with coefficients in $\Z[1/r!]$, and $\P^r$ denote the
subring of $\P$ consisting of power series whose coefficient of $h^d$ is in $P^{d+r}$
for each $d$. Set $\po^r = P^r\cap\po$ and $\ppo^r = \P^r\cap\ppo$.

\pagebreak  

\claim $(j,0)=f(p)$ for some $f$ in $\ppo^{j+3}$.  
\endclaim

Here $f(p)$ denotes the power series obtained by plugging in $p$ for the un-named
variable in each coefficient of $f$.  The desired inequality (15) is immediate from
the claim and the following easy result. 

\lemma{3.8}
If $f\in\ppo^r$ and $f(p)$ has \,{\rm integer} coefficients, then $\osp f(p) \ge
p-r$.
\endlemma

\pf{3.8} Under the hypotheses, the coefficient of $h^d$ in $f(p)$ is an integer
of the form $f_d(p)$ for some {\sl odd} polynomial $f_d$ of degree $\le d+r$ with
coefficients in $\Z[1/(d+r)!]$. If $d<p-r$, then $d+r<p$ and so $f_d(p)$ is
divisible by $p$. Hence $\osp(f(p)) \ge p-r$.
\endpf

To prove the claim, first observe that $[\cdot]\in\P^1$.  In other words the
coefficient $[k]_d$ of $h^d$ in $[k]$ (viewed as a power series) is a polynomial in
$k$ of degree $\le d+1$ with coefficients in $\Z[1/(d+1)!]$. This can be seen using
calculus and the substitution $t=\exp(h/4)$ (see Remark 1.2) or by appealing to the
following useful lemma.  To state it, recall the ``discrete derivative"
$\Delta:P^{r+1}\to P^r$ defined by $\Delta F(\cdot) = F(\cdot + {1\over2}) -
F(\cdot-{1\over2})$.  Now for any polynomial $f \in P^r$, consider the function
$\Sigma f:\Z\to\Q$ given by the sum
$$
\Sigma f(k) = \sum_{j=-m}^{m} f(j)
\tag 16
$$
where $m = (k-1)/2$. The point of the lemma is first to show that discrete
integrals exist and are unique up to constants, and then to use this to show that
$\Sigma f(\cdot)$ is in disguise a polynomial in $\po^{r+1}$.

\lemma{3.9}
Given $f\in P^r$, there exist unique $F,G \in P^{r+1}$ with $G$ odd such that
$\Delta F = f$, $F(0) = 0$, and $G(k) = \Sigma f(k)$ for every integer $k$.
\endlemma

\pf{3.9}
For each integer $d\ge0$ define $F_d\in P^d$ by $F_d(x)=x^d/d$, and set
$f_d=\Delta F_{d+1}$. Clearly $f_d$ is a monic polynomial of degree $d$ with
coefficients in $\Z[1/(d+1)!]$, and so $f = \sum_{d=0}^r a_d f_d$ for suitable
$a_d\in\Z[1/(r+1)!]$. Thus $F=\sum_{d=0}^r a_d F_{d+1}$ is in $P^{r+1}$ with
$\Delta F=f$ and $F(0)=0$. If $E$ is any other polynomial with $\Delta E=f$, then
$\Delta(E-F)=0$. But then the polynomial $E-F$ is periodic, and therefore
constant, whence $E=F$ if $E(0)=0$. 

The sum in (16) defining $\Sigma f$ telescopes when $\Delta F$ is substituted for
$f$, and so equals $F(k/2) - F(-k/2)$. Thus the polynomial $G$ defined by $G(x) =
F(x/2) - F(-x/2)$ satisfies the last equality, and is unique since any polynomial
is determined by its values on the integers.  Clearly $G$ is odd.
\endpf

\pagebreak

The fact that $[\cdot]$ is in $\P^1$ can now be seen as follows:  Lemma 3.9 shows
that $\Sigma$ can be viewed as an operator $P^r\to\po^{r+1}$. This operator extends
coefficient-wise to an operator $\Sigma : \P^r \to
\ppo^{r+1}$.  Substitute $1+h$ for $t$ in $[k] = (t^{2k}-t^{-2k})/(t^2-t^{-2}) =
\sum_{j=-m}^{m} t^{4j}$ and apply the binomial theorem to get
$$
[k] = \sum_{d\ge0} \, \sum_{j=-m}^{m} {4j\choose d} \, h^d
$$ 
(the binomial coefficient is defined for {\sl all} $j$ by
$4j(4j-1)\cdots(4j-d+1)/d!$ if $d>0$, and $1$ if $d=0$).  Evidently
${4\cdot\choose d} \in P^d$ and $[\cdot]_d = \Sigma{4\cdot\choose d}$. Thus
$[\cdot]_d \in \po^{d+1}$, and so $[\cdot] \in \P^1$ by definition.

Now the product of series in $\P$ induces bilinear operators
$\P^r\times\P^s \to \P^{r+s}$, and the $j$th derivative (with respect to
$h$) induces operators $\P^r\to\P^{r+j}$. It follows that the series
$\Sigma([\cdot]^{(j)}[\cdot])$ is in $\ppo^{j+3}$. But $(j,0) =
(\Sigma([\cdot]^{(j)}[\cdot])(p)$, by definition.  This establishes the claim and
thus completes the initial step of the induction.

\subsection{Inductive step}

Differentiating the cabled quantum integers gives $[k,c-1)'=c[2]'[k,c)$, which
leads by a simple calculation to the formula 
$ 
(j,c) = \left((j,c-1)'-(j+1,c-1)\right)/(c[2]'). 
$ 
Now $[2]=t^2+t^{-2}$, so $[2]' = 2t-2t^{-3} = 8h + \cdots$ has $p$-order 1.
Thus $\op(c[2]')=1$, since $\op(c)=0$ (note that $0<c<m$ by assumption).

By the inductive assumption and Lemma 3.4c, $\op(j,c-1)' \ge 2(n-c+1)-j-1 =
2(n-c)-j+1$ and $\op(j+1,c-1) \ge 2(n-c+1)-(j+1) = 2(n-c)-j+1$. Since $\op$ is a
prevaluation, it follows that $\op(j,c) \ge 2(n-c)-j$, proving (14) and thus the
proposition.
\endproof

Combining Propositions 3.5--3.7, we obtain bounds on the orders of the
$p$-brackets of links in the $3$-sphere.

\theorem{3.10}   Let $L$ be a framed link in the $3$-sphere with framings
$a_1,\cdots,a_\ell$.  Then for any odd prime $p=2n+3$,
$$
\op\br L \ge (\ell + {d-1\over d+1}\,b)\,n
$$
where $d$ is the Milnor degree of $L$ and $b$ is the number of framings
divisible by $p$. For $L$ diagonal {\rm($d\ge2$)} and $b>0$, the bound $\op\br
L \ge (\ell+1)n$ holds as well.
\endtheorem

\pagebreak

\proof  Order the components of $L$ so that $p$ divides $a_i$ for $i\le b$. 
For each coloring $c = (c_1,\dots,c_\ell)$ of $L$ set $c_{\max}=\max(c_i)$,
$|c|=\sum_{i=1}^{\ell} c_i$ and $|c|_p = \sum_{i=1}^{b} c_i$.  Since $\op$ is
a prevaluation, it follows from Proposition 1.5 that $\op\br L$ is bounded below
by the minimum over all colorings $0\le c\le n$ of $\sum_{i=1}^{\ell} \op(a_i|c_i)
+ \op(\phi_{L^c})$.  Applying the bounds for $\op(a_i|c_i)$ in Propositions
3.6--3.7 and the bound for $\op(\phi_{L^c})$ in Proposition 3.5 gives
$$
\op\br L \ge (\ell+b)n + \underset{0\le c\le n}\to\min \, 
({d-1\over d+1}|c| - |c|_p).
$$
The minimum is clearly achieved when the first $b$ components have the maximum
allowable cabling index $n$, and the remaining components have index $0$, whence
$|c|=|c|_p=bn$.  The first inequality in the theorem is now immediate.

If $L$ is diagonal, then $\op\br L \ge (\ell+b)n + \min(c_{\max}-|c|_p)$ by
the alternative bound given in Proposition 3.5.  For $b>0$, this expression is
minimized for the same cabling as before, and so  $c_{\max}=n$ and $|c|_p =
bn$.  Therefore in this case $\op\br L \ge (\ell+1)n$.
\endproof

%\vskip -5 pt
\section{Exact values of the $p$-order}

Computational evidence suggests that the bounds for the $p$-orders of the
$p$-sums $\ac$ given in Propositions 3.6 and 3.7 may be sharp:  

\question{Q(a$|$c)}
Let $p=2n+3$ be an odd prime, $a$ and $c$ be integers with $0\le c\le n$. 
Set $r=2$ if $p$ divides $a$, and $r=1$ otherwise.  Is $\op\ac = r(n-c)$? 
\endquestion

Ohtsuki's work \cite{\OI,\S3.6}, together with the improvements in \cite\LW,
suggest a positive answer to $Q(\pm1|c)$; one must show that the invariants
$\nu_{\pm1,c,0}$ of \cite\LW\ are not divisible by $p$, as verified there for
$c\le2$.

The answer to $Q(a|0)$ is yes.  To see this, recall that the $p$-sum $(a|0)$,
or equivalently the $p$-bracket $b_a$ of the $a$-framed unknot, is given by
$$
b_a = \sum_{k=1}^{p/2} t^{a(k^2-1)}[k]^2.
$$ 
Closed forms for these sums as elements of $\lp$ are well known, and can be
obtained by elementary computations involving Gauss sums.  (Up to a factor, $b_a$
coincides with the quantum invariant of the lens space $L(a,1)$; formulas for
general lens spaces in terms of Dedekind sums are known as well \cite\KMII \,
\cite\Je \, \cite\Ga.)  In particular it will be seen below that $b_a$ is an
associate of $h^{rn}$ in $\lp$, and so $\op(b_a)=rn$.  The exact form of $b_a$ will
not concern us here, except when $a\equiv0\pmod p$.

\pagebreak

\proposition{3.11}
Let $p=2n+3$ be an odd prime and $a$ be an arbitrary integer.  Then $\op(b_a) =
rn$, where $r=1$ or $2$ according to whether $a$ is prime to $p$ or not.  In
fact there exists a unit $u_a$ in $\lp$, depending only on the mod $p$ residue
class of $a$, such that 
$
b_a = u_a h^{rn}.
$
Furthermore $(-1)^nu_\oo$ is a square in $\lp$, that is $u_\oo = (-1)^nv_\oo^2$ for
some unit $v_\oo$ in $\lp$.
\endproposition

\proof 
Clearly $[p-k]=-[k]$ and $t^{a((p-k)^2-1)} = t^{a(k^2-1)}$, and so the sum
$\sum_{k=1}^{p/2}$ in the definition of $b_a$ is half of the complete sum
$\sum_{k=1}^p$, denoted simply by $\sum$ below.  

First consider the case when $a$ is prime to $p$, and calculate $b_a$ ``up to
units" (using the notation $x\sim y$ to indicate that $x$ and $y$ are associates). 
Since $h \sim z = s-s^{-1}$, we have $2h^2\,b_a \sim \sum t^{ak^2}(s^{k}-s^{-k})^2
= \sum (t^{ak^2+4k} - 2t^{ak^2} + t^{ak^2-4k})$.  Completing the square gives $\sum
(t^{a(k+2\bar a)^2-4\bar a}-2t^{ak^2}+t^{a(k-2\bar a)^2-4\bar a})$,  where $\bar
a$ is the mod $p$ inverse of $a$.  Since the sum is over a complete set of
residues mod $p$, the quadratic part of each of the three terms in this expression
contributes a Gauss sum $G_a = \sum t^{ak^2}$.  Factoring this out gives
$2G_a(t^{-4\bar a}\!-\!1) \sim 2Gh$,  since $h\sim t^k-1$ for any $k$ prime to $p$
and $G = \pm G_a \sim h^{n+1}$ (see Remark 3.1a).  Thus $b_a \sim G/h
\sim h^n$.  

The exact calculation of $b_\oo$ is easy, involving only the sum of a geometric
series.  We have 
$
2z^2b_\oo = \sum(s^k-s^{-k})^2 = \sum(q^k-2+q^{-k}) = -2p,
$
since the sum of all the $p$th roots of unity vanishes.  By Remark 3.1a, $p =
(-1)^m [m]!^2 z^{2m}$, and so $b_\oo = -p/z^2 = (-1)^n[m]!^2\,z^{2n} = ((-1)^n
[m]!^2 (z/h)^{2n})\, h^{2n}$.  Therefore $u_\oo = (-1)^nv_\oo^2$
for $v_\oo = [m]!(z/h)^n$.
\endproof

\remarks{3.12}
(a) The computation $b_\oo = -p/z^2$ in the preceding proof can be carried out
just as easily for the family of sums 
$$
s_j = \sum_{k=1}^{p/2} [jk][k],
$$ 
for which $s_1 = b_\oo$.  Indeed $2z^2s_j = \sum (s^{(j+1)k} -  s^{(j-1)k} -
s^{-(j-1)k} + s^{-(j+1)k})$.  If $j\not\equiv\pm1 \pmod p$ then $s^{j\pm1}$ are
both primitive $p$th roots of unity, and so the sum vanishes.  For $j\equiv\pm1$,
only two of the four terms survive, which sum to $\mp2p$.  Thus $s_j = \pm b_\oo$
for $j\equiv\pm1 \pmod p$, and $s_j = 0$ otherwise.  These sums will appear in the
examples in \S5. 

\smallskip

(b) Another family of sums which arise in \S5 and are easily computed is defined by
$$
t_a = \sum_{k=1}^{p/2} s^{a(k^2-1)}[k^2].
$$
(Again these depend only on the mod $p$ residue of $a$.)  Proceeding as above, we
have $2zs^at_a = \sum(s^{(a+1)k^2} - s^{(a-1)k^2}) = G_{2(a+1)} - G_{2(a-1)}$,
where $G_0$ is the ``degenerate" Gauss sum $\sum 1^{k^2} = p$.  This difference
is easily analyzed according to the quadratic nature of $a\pm1\pmod p$, and yields
either $\pm G \pm p$ (when $a\equiv\pm1$), $0$ (when $a\pm1$ are either both
quadratic residues or both non-quadratic residues), or $\pm2G$ (otherwise).  In
particular 
$$
t_{\pm1} = {(\pm1)^mG-p \over 2(q^{\pm1}-1)}
$$
where $m=(p-1)/2$.  This has the same $p$-order as $G/h \sim h^n$, and so
$\op(t_{\pm1}) = n$.  Note that $t_1 \ne t_{-1}$, which has an interesting
topological application (see \S5.5).

\smallskip

(c) The $p$-orders of the trivial sums $u = \sum_{k=1}^{p/2}1$ and $v = \sum_{j\le
k=1}^{p/2}1$ are both zero.  Indeed $u = m$ and $v = m(m+1)/2$ are both relatively
prime to $p$.  More generally for any function $f(k)$ or $g(j,k)$, the sums
$u_f = \sum_{k} t^{f(k)} = u + O(h)$ and $v_g = \sum_{j\le k} t^{g(j,k)} = v +
O(h)$ have zero $p$-order.  Note that if $c = \sum_{k} f(k)$ is also prime to $p$,
then $u_f \ne u_{-f}$, since the linear coefficient in $u_{\pm f} = \pm c$.  A
similar statement holds for $v_g$, which has an interesting topological
application (see \S5.5).  
\endremarks

\heading{\S4. Quantum invariants}\endheading

Fix an odd prime $p = 2n+3$ and a closed oriented $3$-manifold $M$.  In this
section we discuss the quantum $SO(3)$ invariant $\tau_p(M)$, which in the
normalization given here is an element of the cyclotomic field $Q_p$, and give a
lower bound for its order in terms of the mod $p$ first betti number of $M$.  We
also introduce a collection of ``finite type invariants" $\tau_p^d$ that dominate
$\tau_p$.  These are studied in more detail in \cite\CMI.

\section{The $3$-manifold invariant}

Let $L$ be an $\ell$-component integrally framed link in the 3-sphere $S^3$.  Write
$\ell_0$, $\ell_+$ and $\ell_-$ for the number of zero, positive and negative
eigenvalues, respectively, of the linking matrix $A_L$ of $L$ (with framings on
the diagonal).  Alternatively $\ell_0$ can be viewed as the nullity of the
integral quadratic form given by $A_L$.  Define the {\it $p$-norm} of $L$ to be the
element
$$
|L| = b_{1}^{\ell_+} b_{-1}^{\ell_-} (b_\oo/h^n)^{\ell_0}
\tag 17
$$ 
in $Q_p$, where $b_a$ is the $p$-bracket $\br{\bigcirc_a}$ of the $a$-framed unknot
and $h=t-1$ as usual (see (10) for the definition of $p$-bracket).  By Proposition
3.11, $|L|$ lies in the ring of integers $\lp$ in $Q_p$ and is of the form
$$
|L| = uh^{n\ell}
\tag 18
$$
for some unit $u\in \lp$ (namely $u_{1}^{\ell_+} u_{-1}^{\ell_-}
u_\oo^{\ell_0}$).  In particular $\op|L| = n\ell$.  

Now consider the 3-manifold $S^3_L$ obtained by surgery on $L$.  Since every
3-manifold arises in this way for some link, we may assume 
$
M = S^3_L.
$ 
Note that the nullity $\ell_0$ is an invariant of $M$, namely the
first Betti number $\b(M) = \rk(H_1(M))$, since $A_L$ is a presentation matrix for
$H_1(M)$.  Similarly the nullity of $A_L$ as a form over $\Z_p$ is the mod $p$ first
Betti number $\bp(M) = \rk(H_1(M;\Z_p))$. If $L$ is diagonal, then $\b(M)$ is the
number of components of $L$ with framing $0$, and $\bp(M)$ is the number with
framings divisible by $p$.    

Define the {\it level $p$ quantum $SO(3)$ invariant} of $M$ to be the quotient
$$
\tau_p(M) = \br L/|L|
\tag 19
$$ 
where $\br L$ is the $p$-bracket viewed as an element of $\lp$.  (Here $p$ need not
be prime, although it must be odd to give a well defined invariant.)  Clearly
$\tau_p(M)$ lies in $\lp[h^{-1}] \subset Q_p$, by (18), and in fact in $\lp$ as
will be seen below.  It is also clear that the invariant $\tau_p$ is {\sl
multiplicative} under connected sums, $\tau_p(M\#N) = \tau_p(M) \tau_p(N)$, since
the bracket is multiplicative and $\ell_0,\ell_\pm$ are additive under distant
unions of links (cf.\ \cite{\KMI,\S5.9}). The motivation for calling it the
``$SO(3)$" invariant arises in the quantum group setting: the $p$-bracket of $L$ can
be rewritten using the ``Symmetry Principle" (\cite{\KMI,\S4.20}) as the sum over
{\sl odd} colorings $k<p$, which correspond to those representations of $SU(2)$ which
factor through $SO(3)$ (cf.\ \cite\BHMV).

\remark{4.1}
This invariant was first defined in \cite{\KMI,\S8.10} with a slightly
different normalization, denoted there by $\tau'_p$.  In particular $\tau'_p(M)
= \br L/|L|'$ (evaluated at $t = \exp(2\pi im^2/p)$) where $|L|'= b_{+1}^{\ell_+}
b_{-1}^{\ell_-} b_{0}^{\ell_0/2}$.  This definition is shown independent of the
choice of $L$ by establishing the invariance of $\br L$ under ``handle-slides",
one of the two moves in the Kirby calculus relating any two framed links which
give the same 3-manifold \cite\Ki.  (See \cite\LII\ for a purely skein theoretic
proof of this invariance.)  It is then an elementary exercise to show that
invariance under ``blow-ups" (the other move) is achieved by dividing by {\sl any}
factor of the form $|L|_c = b_{+1}^{\ell_+} b_{-1}^{\ell_-} c^{\ell_0}$ ($c$ can
be an arbitrary constant since $\ell_0$ is an invariant of $M$).  In particular
the constants for $\tau_p$ and $\tau_p'$ are $c=b_\oo/h^n$ and $c'=b_\oo^{1/2}$.  

Of course the choice of constant affects what properties the quantum invariants
have.  For example $\tau_p'$ is involutive, i.e.\ conjugates under orientation
reversal, whereas $\tau_p$ is not.  On the other hand (and this is what motivated
our choice of normalization) $\tau_p$ takes its values in $Q_p$, whereas $\tau_p'$
does not in general because of the square root in $c'$.  However $\tau_p'$ does
take values in $Q_{4p} = Q_p(\sqrt{-1})$ which differ by units from the
corresponding values of $\tau_p$ (assuming $p$ is prime).  Indeed Proposition 3.11
shows that $c = (-1)^nv_\oo^2h^n$ and $c' = \varepsilon v_\oo h^n$, where
$\varepsilon = (-1)^{n/2}$, so $\tau_p(M) = (\epsilon v_\oo)^{\b(M)}\tau_p'(M)$.

\section{Quantum cyclotomic orders}

Define the {\it quantum $p$-order} $\op(M)$ of $M$ to be $p$-order of
$\tau_p(M)$, that is 
$$
\op(M) = \op(\tau_p(M))
\tag 20
$$
where the $\op$ on the right is the $p$-adic valuation on $Q_p$ defined in
\S3.  Observe that $\op$ is {\sl additive} under connected sums, since
$\tau_p$ is multiplicative.

There is a simple relation between the quantum $p$-order of $M = S^3_L$ and the
$p$-order of the $p$-bracket $\br L$:
$$
\op(M) = \op\br L-n\ell.
\tag 21
$$
Indeed $\op(M) = \op\br L - \op|L|$, since $\op$ is a valuation, and $\op|L| =
n\ell$ by (18).  Now the general bound $\op\br L \ge n\ell$ of Theorem 3.10
implies that $\op(M)\ge0$.  Since $\tau_p(M)\in\lp[h^{-1}]$, this recovers
Murakami's integrality result:

\theorem{4.2} {\rm (\cite\MI\,\cite\MR)}
$\tau_p$ takes values in the cyclotomic ring $\lp$.
\endtheorem

The stronger bounds in Theorem 3.10 lead to corresponding bounds for the
quantum $p$-orders of $3$-manifolds, and then combining these with a result of
Murakami and Ohtsuki, to the universal betti number bound alluded to in the
introduction.  Define the {\it Milnor degree} of a $3$-manifold $M$ to be the {\sl
maximum} Milnor degree of all the (integrally) framed links $L$ for which $M=S^3_L$,
and call
$M$ {\it diagonal} if it has Milnor degree $\ge2$.

\theorem{4.3}  Let $M$ be a closed oriented $3$-manifold and $p=2n+3$ be an odd
prime.  Then
$
\op(M)\ge \b_p(M)\,n(d-1)/(d+1),
$ 
where $\bp(M)=\rk(H_1(M;\Z_p))$ and $d$ is the Milnor degree of $M$.  In
particular for diagonal $M$ {\rm(}meaning $d\ge2${\rm)},
$$
\op(M) \ge \b_p(M)\,n/3.
$$ 
In fact this bound holds for {\rm any} $M$.  Furthermore $\op(M) \ge n$ if
$\b_p(M)$ is positive. 
\endtheorem

\proof
The first inequality is immediate from (21) and Theorem 3.10, as are the other two
for the case when $M$ is {\sl diagonal}.  But by the diagonalizing lemma of
Murakami and Ohtsuki \cite{\MII,\S2.3}, there exists a diagonal $3$-manifold $N$
with $\op(N) = \bp(N) = 0$ such that $M\#N$ is diagonal; in particular $N$ may be
taken to be a connected sum of lens spaces $L(k_i,1)$ with $k_i$ prime to $p$.  Thus
$M$ has the same quantum $p$-order and mod $p$ first betti number as a diagonal
manifold (since both $\op$ and $\bp$ are additive under connected sums) and so both
inequalities hold in general.  
\endproof

\vskip -5 pt
\section{The finite type invariants $\tau_p^d$ and further bounds for the quantum
$p$-order }

Let $\M$ denote the free abelian group generated by closed oriented 3-manifolds
(up to oriented diffeomorphism).  Theorem 4.2 shows that $\tau_p$ can be viewed as
a $\Z$-linear map
$$
\tau_p:\M \to \lp.  
$$
Due to the well known computational complexity of $\tau_p$, one would not expect
this map to be of ``finite type" in any natural sense.  It turns out, however, that
$\tau_p$ is a {\sl limit} of finite type invariants, in the sense of the theory
developed by the authors in \cite\CMI.  

Recall from \cite\CMI\ that a $\Z$-linear map from $\M$ to an abelian group $A$ is
said to be of {\it finite type} if there exists an integer $d\ge0$ such that
$\lambda(M_{\delta L}) = 0$ for all {\it admissible} pairs $(M,L)$ with $\ell>d$;
the smallest such $d$ is called the {\it degree} of $\lambda$.  Here $L$ is an
$\ell$-component framed link in a $3$-manifold $M$, and admissibility means
that each component of $L$ is null-homologous in $M$ with framing $\pm1$ and zero
linking number with any other component.  The notation $M_{\delta L}$ represents the
alternating sum of manifolds obtained by surgery on all the sublinks of $L$, 
$$
M_{\delta L} = \sum_{S<L}(-1)^s M_S
$$
(cf.\ Remark 1.4).  

Now the fact that $\tau_p$ is not of finite type (for $p\ne3$ since $\tau_3\equiv1$)
is an easy consequence of Murakami's beautiful formula  $\tau_p(M) = 1 + 6\lambda(M)h
+ O(h^2)$ for homology spheres, where $\lambda$ is the Casson invariant
\cite\MI.  Indeed it follows from this formula that $\tau_p(S^3_{\delta T_d}) \ne0$
for any $d>0$, where $T_d$  is the distant union of $d$ copies of the $+1$-framed
right-handed trefoil.

It is shown in \cite\CMI, however, that the maps
$$
\tau_p^d:\Cal M \to \Z_{p^k},
$$
obtained from $\tau_p$ by composing with the projections $\pi^d:\lp \to \Z_{p^k}$
defined in (3.2), are of finite type (of degree $\le 3d$, in fact $\le 3d-b(p-3)/2$
when restricted to manifolds with $\bp=b$).  It follows that $\tau_p$ is
``dominated" by finite type invariants.  In other words the value of $\tau_p(M)$ can
be recovered from the values $\tau_p^d(M)$ for all $d$ (just as one can recover an
integer from its residues mod $p^k$ for all $k$).  For $\Zp$-homology spheres $M$,
this can also be deduced from the recent result of Rozansky \cite\RII\ (conjectured by
Lawrence \cite\La) that $\tau_p(M)$ is the $p$-adic limit of the Ohtsuki series
$\sum\lambda_nh^n$
\cite\OI\,\cite\OII.

\remark{4.4}  The degree zero invariants are familiar algebraic topological
invariants.  Indeed the equivalence relation on $3$-manifolds generated by surgery on
admissible links coincides with notion of {\it $H_1$-bordism} \cite\Ge \,
\cite{\CGO,\S2.8,3.1}.  (Recall that $3$-manifolds $M$ and $N$ are $H_1$-bordant if
there is a $4$-manifold $W$ with boundary $M\cup-N$ such that the inclusions $M\to W
\gets N$ induce isomorphisms on $H_1$.)  Thus the degree zero invariants are exactly
the invariants of $H_1$-bordism, including for example the first betti number and
the mod $p$ first betti numbers for each $p$.  In fact the $H_1$-bordism class of a
$3$-manifold has a characterization in terms of its cohomology ring and linking form
(see \cite{\CGO,\S3.1}).  For example, the $H_1$-bordism class of a connected sum of
$b$ copies of $S^1\times S^2$ consists of all $3$-manifolds with $H_1 = \Z^b$ and
with vanishing triple cup product form on $H^1$ \cite{\CGO,\S3.6}.
\endremark

The perspective on the study of $3$-manifolds suggested by the theory of finite type
invariants leads to sharper bounds on the quantum $p$-order.  As an illustration of
this, we prove a result about the orders of manifolds $H_1$-bordant to a connected
sum of $S^1\times S^2$'s.

\proposition{4.5}  Let $M$ be $H_1$-bordant to $\#^b(S^1\times S^2)$.  Then for any
odd prime $p=2n+3$, the quantum $p$-order $\op(M)\ge bn/2$. 
\endproposition

Note that $\#^b(S^1\times S^2)$ has order $bn$ (as will be seen in \S5) while Theorem
4.3 gives a lower bound of $bn/3$ for the order of any manifold with the same first
betti number.  Further more, the examples in the next section realize the lower
bound of $\lceil b/3 \rceil n$ (which is generally much closer to $bn/3$ than to
$bn/2$) for such manifolds.  Thus the bound in the Proposition reflects a strong
restriction on the orders of manifolds $H_1$-bordant (and thus ``closer") to
$\#^b(S^1\times S^2)$, providing further insight on the topological nature of
$\op(M)$.

\pf{4.5}
By the remarks above, $M=S^3_{L}$ with $L=J\cup K$, where $J=\bigcirc^b$ (the
zero-framed unlink of $b$ components) and $K$ is admissible in $S^3_J =
\#^b(S^1\times S^2)$ (i.e.\ $K$ has null homologous components with zero linking
numbers and $\pm1$ framings).  It follows that $K$ is admissible in $S^3$ and that
each component of $K$ has zero linking number with each component of $J$.  Now using
(21), and proceeding as in the proof of Theorem 3.10, we have
$$
\op(M) \ge bn + \underset{0\le c\le n}\to\min \, 
(\o(\phi_{L^c})-|c| - |c|_p). 
$$

As in the proof of Theorem 2.5, the order of $\phi_{L^c}$ is bounded below by the
minimum possible degree of a chordless Feynman diagram on $L^c$ which has (1) at
least two external vertices on each component of $L^c$, and (2) no trees with {\sl
all} their external vertices on the cabling of $J$ (since $J$ has infinite Milnor
degree).  In fact one need only consider diagrams with {\sl exactly} two vertices on
each component of $L^c$, since excess vertices can be traded for loops without
affecting the degree (see Figure 1).  Such a diagram will be called {\it admissible}.

\bigskip
\centerline{\dy{600} \kern 40pt $\longrightarrow$ \kern 40pt \dloop{600}}
\bigskip
\centerline{Figure 1}
\medskip

Set $j=j(c)=2|c|_p$ and $k=k(c)=2(|c|-|c|_p)$, representing the number of external
vertices in admissible diagrams on $L^c$ which lie on the cabling of $J$ and $K$,
respectively.  (Note that $j\le2bn$ since $c\le n$.)  Then
$$
\op(\phi_{L^c}) - |c| - |c|_p \ge k/2 - t(c)
$$ 
where $t(c)$ is the maximum possible number of trees in an admissible diagram on
$L^c$ (cf. the proof of Theorem 2.5).  Since chords are disallowed, $t(c) \le y(c) +
(j+k-3y(c))/4 = (j+k+y(c))/4$ where $y(c)$ is the maximum possible number of $y$'s
(dashed degree 2 trees) in an admissible diagram on $L^c$.  Because of condition (2)
above, every $y$ must have atleast one vertex on the cabling of $K$, and so it is
clear that $y(c)$ is bounded above by the function $f(j,k) = \min(k,(j+k)/3)$.  (Note
that $f$ assumes the first value of the minimum if $j\ge2k$, and the second value if
$j\le2k$.)  Thus $t(c) \le (j+k+f(j,k))/4$.

It follows that
$$
\op(M) \ge bn +  \underset{0\le c\le n}\to\min \, g(j,k)
$$
where $g(j,k) = (k-j-f(j,k))/4$.  One readily computes $g(j,k) = (k-2j)/6$ for
$j\le2k$, and $g(j,k) = -j/4$ for $2k \le j \le 2bn$.  Since $\partial g/\partial j$
is negative, and $\partial g/\partial k$ is positive for $j<2k$ and zero for $j>2k$,
$g$ assumes its minimum value of $-bn/2$ when $j=2bn$ (the maximum allowable value of
$j$) and $k\le bn$.  Therefore $\op(M) \ge bn/2$.
\endproof  
 
\remark{4.6}
The previous proposition can also be deduced from Theorem 4.3, since it is known that
a $3$-manifold $M$ with $H_1(M)$ torsion free is $H_1$-bordant to a connected sum
of $S^1\times S^2$'s if and only the Milnor degree $d(M)\ge 3$ \cite{\CGO, \S\S3.5 and
6.10}.  In a future paper, we will elaborate on this point of view, giving a
characterization $d(M)$ in terms of Massey products.  In particular, it will
be shown that $d(M) = w-1$ where $w$ is the weight of the first nonvanishing Massey
product of $M$ \cite\CMII.
\endremark

\heading{\S5. Examples}\endheading

In this section it is shown that the general betti number bounds established in
Theorem 4.3 are sharp.  To accomplish this, we use examples constructed from three
familiar framed links:  the unknot $L_1$, the (left-handed) Whitehead link $L_2$, and
the Borromean rings $L_3$ (see Figure 2), all equipped with the zero-framing.  

\bigskip
\centerline{\un{500} \kern 30pt \wh{600} \kern 35pt \bo{550} \raise5pt\hbox{\ $=$\ }
\boa{550}}
\bigskip
\centerline{(a) $L_1$ \kern 55pt (b) $L_2$ \kern 90pt (c) $L_3$ \kern 25pt}
\smallskip
\centerline{Figure 2}
\medskip

\noindent 
Observe that $L_2$ and $L_3$ have Milnor degree $3$ and $2$, respectively, while
$L_1$ has infinite Milnor degree.  It follows from Theorem 4.3 that the quantum
$p$-orders of the $3$-manifolds $M_\ell = S^3_{L_\ell}$, obtained by surgery on
these framed links, are all $\ge n$.  In fact these manifolds are all of $p$-order
{\sl exactly $n$}, as will be seen below.  (This result is well known for
$M_1=S^1\times S^2$ and $M_3 = T^3$.)

\theorem{5.1}
For any odd prime $p=2n+3$, the $3$-manifolds $M_1$, $M_2$ and $M_3$ all have
quantum $p$-order $n$.   
\endtheorem

Since $\bp(M_\ell) = \ell$, it follows that the bound $\op\ge n$ is sharp for $\bp =
1$, $2$ or $3$.  Furthermore $\op(\#^kM_3) = kn$, since $\op$ and $\bp$ add
under connected sums.  Therefore the bound $\op\ge \bp n/3$ is sharp for
$\bp\equiv0\pmod3$.

\remark{5.2} It is not known whether the betti number bound $\op\ge\bp n/3$ is sharp
for all $\bp$.  To show this, it would suffice to produce two $3$-manifolds $M_\ell$
(for $\ell=4,5$) with $\bp(M_\ell) = \ell$ and $\op(M_\ell) = \lceil \ell n/3
\rceil$, where $\lceil \, \rceil$ is the least integer upper bound function; 
possible candidates are surgeries on zero framed links $L_\ell$ with $\bar \mu_{ijk}
= 1$ for all distinct triples $i,j,k$ (i.e.\ every three component sublink has a
``Borromean interaction").  The sharpest betti number bound that can be deduced
from the examples in this paper is $o_p \ge \lceil b/3 \rceil n$, realized by
connected sums of copies of $M_2$ (at most two) and $M_3$.
\endremark

\pf{5.1}  By (21) it suffices to show $\op\br{L_\ell} = (\ell+1)n$ for $\ell = 1$,
$2$ and $3$.  This is immediate for $\ell = 1$ since $\br{L_1} = b_\oo$, which has
$p$-order $2n$ by Proposition 3.11.  For $\ell=2$ and $3$, this will be
accomplished by expressing $\br{L_\ell} \in \lp$ in terms of $b_\oo$ and the sums
discussed in Remark 3.12.  In particular, it will be shown that $\br{L_2} = b_\oo
t_1$ and $\br{L_3} = b_\oo^2 u$, which will prove the theorem.

Recall from (10) that the $p$-bracket $\br L$ of any zero-framed link $L$ is given
by the linear combination $\sum_{k=1}^{p/2}[k]J_{L,k}$ of colored Jones polynomials
$J_{L,k}$.  By allowing link colorings in the group ring $\lp\Z$ and expanding
multilinearly,  $\br L$ can be viewed as a single colored Jones polynomial
$J_{L,\omega}$ with each component colored by 
$$
\omega = \sum_{k=1}^{p/2} [k]k \, \in \, \lp\Z.
$$ 
This point of view was introduced by Lickorish \cite\LI\ and will be
adopted here.  In particular a {\it coloring} will henceforth mean a
$\lp\Z$-coloring (which of course need not assign the same color to each
component).

Now consider the equivalence relation $\approx$ on the set of $\lp$-linear
combinations of colored links, defined by $(L,\lambda) \approx (L',\lambda')$ iff
$J_{L,\lambda} = J_{L',\lambda'}$.  For example the following local equivalences
involving integer colors (indicated by labels $j$ and $k$) are well known:
$$
\text{a)} \quad \unt{200}
\put{20}{10}{$\ssize k$}
\ \strand{j} 
\quad \approx \quad 
[k]\ \strand{j}
\kern 60pt
\text{b)} \quad \uni{300}
\put{26}{9}{$\ssize k$}
\put{16.5}{0}{\strandi{j}} 
\quad \approx \quad 
{[jk]\over [j]}\ \strand{j} 
\tag 22
$$
\smallskip\noindent
cf.\ \cite{\KMI,\S3.27}.  For simplicity it is assumed that $0<j,k<p/2$.

It follows that
$$
\text{a)} \quad \unt{200}
\put{20}{10}{$\ssize\omega$} 
\ \strand{j} 
\quad \approx \quad 
b_\oo \ \strand{j}
\kern 60pt
\text{b)} \quad \uni{300}
\put{26}{9}{$\ssize\omega$}
\put{16.5}{0}{\strandi{j}} 
\quad \approx \quad 
\delta_{1j} \ b_\oo
\tag 23
$$
\smallskip\noindent
where $\delta$ is the Kronecker delta and $b_\oo$ is the $p$-bracket of the
zero-framed unknot.  Indeed the (23a) is immediate from (22a) and the definition of
$\omega$.  The equivalence (23b) says that the colored Jones polynomial of the
left-hand link vanishes unless $j=1$, in which case it equals $b_\oo$ times the
polynomial of the link obtained by removing the two pictured components.  This
follows from Remark 3.12a and the obvious fact that $1$-colored components can be
ignored, since the left hand side reduces to $[j]^{-1}s_j\,|^j$ using (22b) (also
see \cite{\LII, Lemma 6}).  

The following generalization of (23) provides the key to calculating
the $p$-bracket $\br{L_\ell} = J_{L_\ell,\omega}$:
$$
\unii{300}
\put{35}{8}{$\ssize\omega$}
\put{27.5}{0}{\strandi{i}}
\put{11.5}{0}{\strandi{j}} 
\quad \approx \quad 
{\delta_{ij}\over(2w,j]} \,b_\oo\ \caps{300}
\put{12}{17}{$\ssize j$}
\put{12}{-13}{$\ssize j$}
\tag 24
$$
\smallskip\noindent
Here $(2w,j]$ is the framed quantum integer $s^{w(j^2-1)}[j]$, where $w$ is defined
by the following scheme.  Observe that in the left-hand link, the two vertical
strands are either a) oppositely oriented arcs on the same component, b)
identically oriented arcs on the same component, or c) arcs from distinct
components.  For each of these cases in turn, define $w$ to be
$$
\text{a)} \quad \lk(\,\orcaps{150}\,)
\kern 50pt
\text{b)} \quad 2\,\lk(\,\cross{150}\,)-1
\kern 50pt
\text{c)} \quad \lk(\,\down{150}\,)
$$  
where the picture specifies a two component oriented link formed by locally
modifying the component(s) containing the vertical strands.

The framed version of (24) is well known and has appeared in various forms in the
literature (see for example \cite{\Tu,\S II.3.10}).  In particular, one defines
the Jones polynomial of a {\sl framed} colored link $(L,k)$ to be $t^{a(k^2-1)}
J_{L,k}$, where $a$ is the framing, with associated equivalence relation
$\approxf$.  Then (24) takes the diagramatic form 
$$
\unii{300}
\put{36}{8}{$\ssize\omega$}
\put{27.5}{0}{\strandi{i}}
\put{11.5}{0}{\strandi{j}} 
\quad \approxf \quad 
{\delta_{ij}\over[j]} \,b_\oo\ \caps{300}
\put{12}{17}{$\ssize j$}
\put{12}{-13}{$\ssize j$}
\tag 25
$$
\smallskip\noindent
with respect to the ``blackboard framing", since the sum of the blackboard
framings of the components on the left-hand side exceeds the corresponding sum on
the right-hand side by exactly $2w$.  

The proof of (25) is easiest in the general context of framed colored {\sl
tri-valent graphs} (cf.\ \cite\KR \, \cite\KL).  In the quantum group approach, the
Jones polynomial of such a graph is defined locally in terms of operators assigned
to the elementary tangles $|$, $\cup$, $\cap$, $\pos{100}$, $\neg{100}$, $\inj{100}$
and $\pro{100}$.  Specifically $\inj{100}$ represents a natural injection of a
simple $U_q(sl_2)$-module into the tensor product of two other simple modules,
and $\pro{100}$ represents the corresponding projection.  These satisfy the
following identities:
$$
\text{a)} \quad 
\strand{i} \kern 12pt \strand{j} 
\quad \approxf \quad  
\sum_m \ \spi{300}
\put{21}{19}{$\ssize i$} \put{20}{-18}{$\ssize i$} \put{1}{19}{$\ssize j$}
\put{1}{-18}{$\ssize j$} \put{7}{0}{$\ssize m$}
\kern 50pt
\text{b)} \quad
\spi{300} 
\put{22}{19}{$\ssize j$} \put{22}{-18}{$\ssize j$} 
\put{1}{19}{$\ssize j$} \put{1}{-18}{$\ssize j$} 
\put{7}{0}{$\ssize 1$} 
\quad \approxf \quad  
\dsize{1\over[j]} \ \caps{300}
\put{12}{17}{$\ssize j$}
\put{12}{-13}{$\ssize j$}
\tag 26
$$
\smallskip \noindent
where the sum is over all {\it admissible} $m$, i.e.\ $|i-j| < m < i+j$ with $i+j+m$
odd (see \cite{\KR,\S4.11-12}).  By (26a)
$$
\unii{300}
\put{36}{8}{$\ssize\omega$}
\put{27.5}{0}{\strandi{i}}
\put{11.5}{0}{\strandi{j}} 
\quad \approxf \quad 
\sum_m 
\ \ \spis{300} \put{31}{7}{$\ssize\omega$}
\put{26}{19}{$\ssize i$} \put{26}{-18}{$\ssize i$} 
\put{2}{19}{$\ssize j$} \put{2}{-18}{$\ssize j$} 
\put{12}{0}{$\ssize m$} 
$$
\smallskip\noindent
which vanishes by (23b) unless $m=1$.  But $1$ is admissible if and only if $i=j$,
and so (25) follows from (23b) and (26b).

It is now straightforward to compute the $p$-brackets of $L_2$ and $L_3$.  Applying
(24) to $(L_2,\omega)$ gives   
$$
^\omega\wh{400}
\put{0}{10}{$\ssize\omega$} 
\quad = \quad
\sum_k \, [k] \ ^\omega\wh{400}
\put{0}{10}{$\ssize k$} 
\quad \approx \quad
b_\oo \, \sum_k s^{k^2-1}  \ \hop{400} 
\put{0}{10}{$\ssize k$} \put{0}{-10}{$\ssize k$}
$$
\smallskip \noindent
But $J_{^k\hopf{80}^k} = [k^2]$ by (22), and so $\br{L_2} = b_\oo t_1$.  Similarly,
applying (24) twice to $(L_3,\omega)$ gives
$$
^\omega\boa{400} \, {\ssize\omega} \put{20}{20}{$\ssize\omega$}
\ \ = \ \
\sum_k \, [k] \ ^\omega\boa{400} \, {\ssize\omega}
\put{20}{20}{$\ssize k$}
\ \  \approx \ \ 
b_\oo \, \sum_k \, \tri{400}
\put{2}{15}{$\ssize k$} \put{2}{-15}{$\ssize k$} \, {\ssize\omega}
\ \ \approx \ \  
b_\oo \, \sum_k \, {1\over[k]} \, \bigcirc^k
$$
\smallskip \noindent
and so $\br{L_3} = b_\oo^2u$ \, since $J_{\bigcirc^k} = [k]$.
\endproof

\smallskip
Given a $3$-manifold $M$, let $\rp(M)$ denote the {\sl minimum} quantum $p$-order of
any manifold in its {\it betti class} (all manifolds with the same first betti number
as
$M$).  The previous theorem gives a partial determination of $\rp(M)$, namely $\rp(M)
= \b(M)n/3$ if $\b(M)$ is divisible by $3$, and $\rp(M)=n$ if $\b(M)\le3$, where as
usual $n=(p-3)/2$.   

\remark{5.3}
There is no {\sl maximum} quantum $p$-order for manifolds in a given betti class. 
Indeed for any prime $p$, there exist manifolds $N_b$ with first betti number $b$ and
infinite $p$-order (i.e.\ $\tau_p(N_b) = 0$).  For example, if $p\equiv1\pmod3$ then
let $N_1$ be $0$-surgery on the (right-handed) trefoil, or equivalently $(1,0)$-surgery
on the (right-handed) Whitehead link.  The $p$-bracket of this framed link can be
calculated, as for $\br{L_2}$ in the proof above, to be $b_\oo t_{(p-1)/2}$, which has
infinite $p$-order by Remark 3.12b.  (Note that $(p-1)(p+3)$ is a square mod $p$
under the congruence assumption on $p$, by an elementary calculation using quadratic
reciprocity.)  Now let $N_b = \#^bN_1$ for $b>0$, and $N_0 = p$-surgery on the
trefoil.     
\endremark

It is also an interesting problem to determine the minimum quantum $p$-order of any
manifold in the {\it $H_1$-bordism class} (see Remark 4.4) of $M$, denoted $\sp(M)$. 
For example consider the family of links $L_3^k$ obtained from the Borromean rings by
replacing one of the components with its $(k,1)$-cable ($L_3^{-2}$ is pictured in
Figure 3a) and let $M_3^k$ denote the $3$-manifold obtained by zero-framed surgery on
$L_3^k$.  It is known that $M_3^j$ and $M_3^k$ are $H_1$-bordant if and only if
$j=\pm k$, and that these manifolds represent all the $H_1$-bordism classes of
manifolds with first homology $\Z^3$ (\cite{\CGO,\S3.3}).  By Theorem 4.3 and
Proposition 4.5, $\sp(M_3^k)\ge n$ for all $k\ne0$ and $\sp(M_3^0)\ge 3n/2$.  We
suspect that all of these bounds are sharp.

%\pagebreak

\bigskip
\centerline{\dboat{700} \kern 75pt \ring{300}}
\bigskip
\centerline{(a) \kern 110pt (b)  }
\smallskip
\centerline{Figure 3}
\medskip

Of course $M_3^1$ is just the $3$-torus, and so $\sp(M_3^1)=n$ by the previous
theorem.  In contrast, $M_3^0=\#^3S^1\times S^2$ has order $3n$, which is not
minimal in its $H_1$-bordism class.  Indeed the manifold $S^3_{L_1\sqcup L_2} = M_1\#
M_2$ has order $2n$.  Presumably the $3$-manifold obtained by zero surgery on the
link shown in Figure 3b has $p$-order $\lceil 3n/2 \rceil$, although this has not
been verified.  This would show that $\sp(M_3^0) = \lceil 3n/2 \rceil$.  For
$|k|>1$, we suspect that
$M_3^k$ has $p$-order $n$, which would give $\sp(M_3^k) = n$ for $k\ne0$.  This can
be verified for $k=\pm2$ as follows.

\proposition{5.4} For any odd prime $p=2n+3$, the manifolds $M_3^{k}$ have quantum
$p$-order $n$ for $k=\pm2$.
\endproposition

\proof By (21) it suffices to show $\br{L_3^{k}} = 4n$.  Proceeding as in the proof
of Theorem 5.1, this bracket (for $k=-2$) is equal to the colored Jones polynomial of
$$
^\omega\dboa{400} \, {\ssize\omega} \put{20}{20}{$\ssize\omega$}
\ \ = \ \
\sum_k \, [k] \ ^\omega\dboa{400} \, {\ssize\omega}
\put{20}{20}{$\ssize k$}
\ \  \approx \ \ 
b_\oo \, \sum_k \, \dtri{400}
\put{4}{15}{$\ssize k$} \put{4}{-15}{$\ssize k$} \, {\ssize\omega}
\ \ \approx \ \  
b_\oo \, \sum_{j\le k} \, ^{2j-1} \dhopf{400} \, ^{\omega}
$$ 
\smallskip\noindent
where the final sum is over $1\le j\le k < p/2$.  The last equivalence
follows from a well known cabling principle for colored Jones polynomials (see for
example \cite{\KMI,\S3.10}) and the decomposition of tensor products of simple
$U_q$-modules (\cite{\KMI,\S2.13}).  Since the final pictured link is symmetric, the
equivalence (24) can be applied once more to give $\br{L_3^{\pm2}} = b_\oo^2
\sum_{j\le k} q^{\mp 2j(j-1)}$, which has $p$-order $4n$ by Remark 3.12c.
\endproof

\remark{5.5} The two $3$-manifolds $M_3^2$ and $M_3^{-2}$ (which are orientation
reversing diffeomorphic) not only have the same $p$-orders, but also the same lowest
order coefficients in their quantum invariants.  They are however distinguished by
the next highest order coefficient, for $p\ne3$, and therefore by the finite type
invariant $\tau_p^{n+1}$.  (To see this, one can work with the renormalized invariant
$u_\oo\tau_p/h^n$, where $u_\oo$ is the unit $b_\oo/h^{2n}$ of Proposition 3.11,
which assumes the value $\sum_{j\le k} q^{\mp 2j(j-1)}$ on $M_3^{\pm2}$.  The
constant coefficients in these sums are both equal to $\sum_{j\le k}1 = m(m+1)/2$,
while the linear coefficients $\sum_{j\le k} \mp 8j(j-1) = \mp(m-1)m(m+1)(m+2)/3$
are distinct, since they are negatives and prime to $p$; see Remark 3.12c.) 
Therefore $M_3^2$ is {\it chiral}, that is has no orientation reversing
automorphisms.  A similar argument using Remark 3.12b shows that the manifold $M_2$,
obtained by zero surgery on the Whitehead link, is chiral (also proved easily using
Lescop's generalization of the Casson invariant).
\endremark

\bigskip

\Refs
\def\rf{\smallskip\ref\key}

\widestnumber\key{BHMV}

\rf BHMV
\by C. Blanchet, N. Habegger, G. Masbaum and P. Vogel 
\paper Remarks on the three-manifold invariants $\theta_p$, 
\nofrills {\rm in}
\inbook ``Operator Algebras, Mathematical Physics, and Low Dimensional Topology"
\eds R. Herman and B. Tanbay 
\bookinfo Res. Notes Math.
\vol 5
\yr 1993
\pages 39--59
\endref

\rf BS
\by Z. Borevich and I. Shafarevich 
\book Number Theory
\publ Academic Press 
\publaddr New York
\yr 1966
\endref

\rf Co
\by T.D. Cochran 
\paper Link concordance invariants and homotopy theory
\jour Invent. Math.
\vol 90 
\yr 1987
\pages 635--645
\endref

\rf CGO
\by T.D. Cochran, A. Gerges and K. Orr 
\paper Dehn surgery equivalence relations on three-manifolds 
\jour Rice Univ. preprint  
\vol  
\yr 1997
\pages 
\endref

\rf CM1
\by T.D. Cochran and P. Melvin 
\paper Finite type invariants of $3$-manifolds 
\jour 
\vol  
\yr 
\pages 
\toappear
\finalinfo math.GT/9805026
\endref

\rf CM2
\bysame  % T.D. Cochran and P. Melvin 
\paper The Milnor degree of a $3$-manifold 
\jour   
\vol  
\yr
\pages 
\finalinfo (in preparation)
\endref

\rf Dr 
\by V.G. Drinfel'd 
\paper Quasi-Hopf algebras 
\jour Leningrad Math. J.
\vol 1
\yr 1990
\pages 1419--1457
\endref

\rf Ga
\by S. Garoufalidis
\paper Relations among $3$-manifold invariants 
\jour Thesis, Univ. Chicago
\yr 1992
\endref

\rf Ge
\by A. Gerges
\paper Surgery, bordism and equivalence of $3$-manifolds 
\jour Thesis, Rice Univ.
\yr 1996
\endref

\rf HM 
\by N. Habegger and G. Masbaum 
\paper The Kontsevich integral and Milnor's invariants
\toappear
\jour Topology
\vol  
\yr 
\pages 
\endref

\rf IR
\by K. Ireland and M. Rosen 
\book A Classical Introduction to Modern Number Theory
\publ Springer-Verlag 
\publaddr New York
\yr 1990
\endref

\rf Je
\by L. Jeffrey
\paper Chern-Simons-Witten invariants of lens spaces and torus bundles, and the
semi-classical approximation 
\jour Commun. Math. Phys.
\vol 147
\yr 1992
\pages 563--604
\endref

\rf Jo 
\by V.F.R. Jones
\paper A polynomial invariant of knots via von Neumann algebras
\jour Bull. Amer. Math. Soc.
\vol 12
\yr 1985
\pages 103--112
\endref

\rf Ka 
\by L.H. Kauffman
\paper State models and the Jones polynomial
\jour Topology
\vol 26
\yr 1987
\pages 395--407
\endref

\rf KL 
\by L.H. Kauffman and S.L. Lins
\book Temperley-Lieb Recoupling Theory and Invariants of $3$-Manifolds
\bookinfo Annals of Math. Studies 134 
\publ Princeton Univ. Press
\publaddr Princeton, NJ
\yr 1994
\endref

\rf Ki 
\by R. Kirby
\paper A calculus for framed links
\jour Invent. Math.
\vol 45
\yr 1978
\pages 35--56
\endref

\rf KM1 
\manyby R. Kirby and P. Melvin 
\paper The $3$-manifold invariants of Witten and Reshetikhin-Turaev for $sl(2,\C)$
\jour Invent. Math.
\vol 105
\yr 1991
\pages 473--545
\endref

\rf KM2 
\bysame % R. Kirby and P. Melvin 
\paper Quantum invariants of lens spaces and a Dehn surgery formula
\jour Abstracts Amer. Math. Soc.
\vol 12
\yr 1991
\page 435
\endref

\rf KR 
\by A.N. Kirillov and N.Yu. Reshetikhin
\paper Representations of the algebra $U_q(sl_2)$, $q$-orthogonal polynomials and
invariants of links, 
\nofrills {\rm in}
\inbook ``Infinite Dimensional Lie Algebras and Groups"
\ed V. G. Kac
\publ World Scientific
\publaddr Singapore
\yr 1988
\pages 285--339
\endref

\rf KS
\by A. Kricker and B. Spence 
\paper Ohtsuki's invariants are of finite type
\jour Queen Mary and Westfield Coll. preprint
\vol  
\yr 1996
\pages 
\finalinfo math.QA/9608007
\endref

\rf La
\by R. Lawrence
\paper Asymtotic expansions of Witten-Reshetikhin-Turaev invariants for some simple
$3$-manifolds
\jour J. Math. Phys.
\vol 36
\yr 1995
\pages 473--545
\endref

\rf Le
\by T.T.Q. Le 
\paper An invariant of integral homology $3$-spheres which is universal for all
finite type invariants
\jour Amer. Math. Soc. Transl. (2)
\vol 179
\yr 1997 
\pages 75--100
\endref

\rf LM1
\manyby T.T.Q. Le and J. Murakami  
\paper Kontsevich's integral for the homfly polynomial and relations between
values of multiple zeta functions
\jour Topology Appl. 
\vol 62
\yr 1995 
\pages 193--206
\endref

\rf LM2
\bysame % T.T.Q. Le and J. Murakami
\paper Parallel version of the universal Vassiliev-Kontsevich invariant
\jour J. Pure and Appl. Alg.
\vol 121
\yr 1997
\pages 271--291
\endref

\rf LMO
\by T.T.Q. Le, J. Murakami and T. Ohtsuki
\paper On a universal perturbative invariant of $3$-manifolds
\jour Topology
\vol 37
\yr 1998
\pages 539--574
\endref

\rf L1
\by W.B.R. Lickorish
\paper Calculations with the Temperley-Leib Algebra
\jour Commun. Math. Helv.
\vol 67
\yr 1992
\pages 571--591
\endref

\rf L2
\bysame % W.B.R. Lickorish  
\paper The skein method for three-manifold invariants
\jour J. Knot Theory Rami.
\vol 2
\yr 1993 
\pages 171--194
\endref

\rf LW
\by X.-S. Lin and Z. Wang 
\paper On Ohtsuki's invariants of integral homology $3$-spheres 
\toappear
\jour Acta Math. Sinica
\vol
\yr 
\pages
\endref

\rf MR
\by G. Masbaum and J.D. Roberts
\paper A simple proof of integrality of quantum invariants at prime roots of unity
\jour Math. Proc. Cam. Phil. Soc.
\vol 121
\yr 1997
\pages 443-454
\endref

\rf MW
\by G. Masbaum and H. Wenzl
\paper Integral modular categories and integrality of quantum invariants at roots of
unity of prime order
\toappear
\jour Crelle's Jour. 
\vol
\yr 
\pages
\endref

\rf Mi
\by J. Milnor
\paper Isotopy of links, 
\nofrills {\rm in}
\inbook ``Algebraic Geometry and Topology.  A Symposium in Honor of S. Lefshetz"
\eds Fox, Spencer and Tucker
\publ Princeton Univ. Press
\publaddr Princeton, NJ
\pages 280--306
\yr 1957
\endref

\rf M1
\manyby H. Murakami 
\paper Quantum $SU(2)$-invariants dominate Casson's $SU(2)$-invariant
\jour Math. Proc. Camb. Phil. Soc.
\vol 115
\yr 1994
\pages 253--281
\endref

\rf M2
\bysame % H. Murakami 
\paper Quantum $SO(3)$-invariants dominate the $SU(2)$-invariant of Casson and Walker
\jour Math. Proc. Camb. Phil. Soc.
\vol 117
\yr 1995
\pages 237--249
\endref

\rf O1
\manyby T. Ohtsuki
\paper A polynomial invariant of integral homology $3$-spheres
\jour Math. Proc. Camb. Phil. Soc.
\vol 117
\yr 1995
\pages 83--112
\endref

\rf O2
\bysame % T. Ohtsuki
\paper A polynomial invariant of rational homology $3$-spheres
\jour Invent. Math.
\vol 123
\yr 1996
\pages 241--257
\endref

\rf RT1
\manyby N.Yu. Reshetikhin and V.G. Turaev
\paper Ribbon graphs and their invariants derived from quantum groups
\jour Commun. Math. Phys.
\vol 127
\yr 1990
\pages 1--26
\endref

\rf RT2
\bysame % N.Yu. Reshetikhin and V.G. Turaev
\paper Invariants of $3$-manifolds via link polynomials and quantum groups
\jour Invent. Math.
\vol 103
\yr 1991
\pages 547--597
\endref

\rf R1
\manyby L. Rozansky
\paper Reshetikhin's formula for the Jones polynomial of a link: Feynman diagrams and
Milnor's linking numbers
\jour J. Math. Phys.
\vol 35
\yr 1994
\pages 5219--5246
\endref

\rf R2
\bysame % L. Rozansky
\paper On $p$-adic properties of the Witten-Reshetikhin-Turaev invariant
\jour Yale Univ. preprint
\vol 
\yr 1998
\pages 
\finalinfo math.QA/9806075
\endref

\rf TY
\by T. Takata and Y. Yokota
\paper The $PSU(N)$ invariants of $3$-manifolds are polynomials
\jour Kyushu Univ. preprint
\vol
\yr 1997
\pages
\endref

\rf Tu
\by V.G. Turaev
\book Quantum Invariants of Knots and $3$-Manifolds
\bookinfo Studies in Math.
\vol 18
\publ Walter de Gruyter
\publaddr Berlin
\yr 1994
\endref

\rf Wa
\by L.C. Washington
\book Introduction to Cyclotomic Fields,
\publ Springer-Verlag
\publaddr New York
\yr 1982
\endref

\rf Wi
\by E. Witten
\paper Quantum field theory and the Jones polynomial
\jour Commun. Math. Phys.
\vol 121
\yr 1989
\pages 351-399
\endref

\bigskip\bigskip\sl

\noindent
Department of Mathematics

\noindent
Rice University

\noindent
Houston, TX 77005-1892

\noindent\rm
cochran\@math.rice.edu

\bigskip\sl

\noindent
Department of Mathematics

\noindent
Bryn Mawr College

\noindent
Bryn Mawr, PA 19010-2899

\noindent\rm
melvin\@brynmawr.edu

\end